\input amstex
\documentstyle{pspum}

\loadeusm

\input epsf.tex

\def\lam{{\lambda}}
\def\Lam{{\Lambda}}

 \def\bfg{\bold{g}} 
  
 \def\bft{\bold{t}} 
 
\def\bfp{\bold{p}}

\def\bfg{\frak{g}}

\def\abar{\bar{a}}

\def\barz{\bar{z}}
\def\barsig{\bar{\sigma}}
\def\barlam{\bar{\lambda}}
\def\var{\varphi}

\def\dbA{{\Bbb A}}
\def\dbC{{\Bbb C}}
\def\dbH{{\Bbb H}}
\def\dbP{{\Bbb P}}
\def\dbR{{\Bbb R}}
\def\dbZ{{\Bbb Z}}

\def\scr#1{{\fam\eusmfam\relax#1}}

   \def\scrF{{\scr F}}

   \def\scrL{{\scr L}}
\def\scrM{{\scr M}}   
\def\scrO{{\scr O}}

   \def\scrX{{\scr X}}

\def\gr#1{{\fam\eufmfam\relax#1}}


	\def\grg{{\gr g}}

\def\grP{{\gr P}}

\def\SU{\text{\rm SU}} \def\QFT{\text{\rm QFT}}
\def\SO{\text{\rm SO}} \def\ad{\text{\rm ad}}
\def\SYM{\text{\rm SYM}} \def\Tr{\text{\rm Tr}}
\def\Sym{\text{\rm Sym}} \def\eff{\text{\rm eff}}
\def\SUSY{\text{\rm SUSY}}  \def\susy{\text{\rm SUSY}}
\def\im{\text{\rm Im}} \def\SL{\text{\rm SL}}
\def\Sp{\text{\rm Sp}} \def\sing{\text{\rm sing}}
\def\Loop{\text{\rm loop}} \def\cl{\text{c$\ell$}}

\def\noi{\noindent}
\def\op{\operatornamewithlimits}

\def\dashright{\hbox{-\ -\ -\ -$\kern-3pt$
     \hbox{\lower.1pt\hbox{$\succ$}}}}
\def\w{{\mathchoice{\,{\scriptstyle\wedge}\,}
  {{\scriptstyle\wedge}}
  {{\scriptscriptstyle\wedge}}{{\scriptscriptstyle\wedge}}}}
\def\Le{{\mathchoice{\,{\scriptstyle\le}\,}
  {\,{\scriptstyle\le}\,}
  {\,{\scriptscriptstyle\le}\,}{\,{\scriptscriptstyle\le}\,}}}

\def\upvee{\sssize{\vee}}

\topmatter
\title
Seiberg-Witten Integrable Systems
\endtitle
\author
Ron Y. Donagi
\endauthor
\address
Department of Mathemtics,\hfill\break
\indent University of Pennsylvania,\hfill\break
\indent Philadelphia, PA \ 19104-6395
\medskip
\endaddress
\curraddr
School of Mathematics,\hfill\break
\null\kern3.0cm Institute for Advanced Study,\hfill\break
\null\kern3.0cm Princeton, NJ \ 08540
\endcurraddr
\thanks
 \endgraf
Secondary: 14D20, 14H40. \endgraf
Work partially supported by NSF grant
DMS95-03249, and by grants from the \endgraf
Univesity of Pennsylvania Research
Foundation and The Harmon Duncombe Foundation.
\endthanks
\toc
\widestnumber\head{\S{44}}
\widestnumber\subhead{XX}
\specialhead{} Introduction\endspecialhead
\medskip
\head \S1. Some physics background\endhead
\subhead{} Electromagnetism\endsubhead
\subhead{} Yang-Mills theory\endsubhead
\subhead{} Quantization\endsubhead
\subhead{} Supersymmetry\endsubhead
\subhead{} $N=2$ super Yang-Mills\endsubhead
\subhead{} Duality\endsubhead
\subhead{} The Seiberg-Witten solution\endsubhead
\subhead{} Adding Matter\endsubhead
\bigskip
\head \S2. Why integrable systems?\endhead
\subhead{} Algebraically integrable systems\endsubhead
\subhead{} Seiberg-Witten differentials\endsubhead
\subhead{} Linearity: complexified Duistermaat-Heckman\endsubhead
\bigskip
\head \S3. Which integrable system?\endhead
\subhead{} Meromorphic Higgs bundles\endsubhead
\subhead{} The spectral curves\endsubhead
\subhead{} Elliptic solitons\endsubhead
\subhead{} Tests\endsubhead
\subhead{} Consistency with Seiberg-Witten\endsubhead
\subhead{} Mass to zero: the N=4 limit\endsubhead
\subhead{} $\tau$ to $\infty$: the flow to pure N=2\endsubhead
\subhead{} Higgs to $\infty$: symmetry breaking\endsubhead
\subhead{} Singularities\endsubhead
\bigskip
\head \S4. Other integrable systems\endhead
\subhead{} Pure $N=2$ SYM\endsubhead
\subhead{} The Toda system\endsubhead
\subhead{} Periodic Toda\endsubhead
\subhead{} The spectral curves\endsubhead
\subhead{} Fundamental matter\endsubhead
\subhead{} Adjoint matter\endsubhead
\subhead{} Stringy integrable systems\endsubhead
\bigskip
\specialhead{} References\endspecialhead
\endtoc
\endtopmatter

\document

\head
Introduction
\endhead

While many mathematicians are aware of the monopole equations 
\cite{W}  and of the impact they have had in topology, probably 
not very many are as familiar with the works  \cite{SW1, SW2} 
of Seiberg and Witten on which the monopole equations are based. 
The appearance of these works was a major breakthrough in Quantum 
Field Theory (QFT), and has to a large extent revolutionized and 
transformed the subject. The actual achievement in \cite{SW1} 
was an {\it exact}
solution to a problem in QFT (pure $\ssize N=2$ supersymmetric
Yang-Mills in four dimensions with gauge group $SU(2)$)
which could traditionally be solved only approximately, by
computing the first few terms in a perturbative expansion.
The solution is close to algebraic geometry both in its
details (involving families of elliptic curves) and in its
general approach, which is to use global properties
(monodromy, symmetries) of a physical theory to uniquely
characterize a geometric model of the theory. 
This work opened the doors to solutions, along similar
lines, of many specific quantum theories, some of which are
reviewed here.
It also led directly to the ongoing revolution in string
theory involving string dualities and the appearance of
higher dimensional objects (branes) on more or less equal
footing with the strings themselves.

The purpose of these notes is to try to make this exciting
and rapidly evolving chapter of high energy physics
accessible to mathematicians, particularly to algebraic
geometers.
I have had to make some choices: I decided to focus on the
original work and its further development within QFT, thus
omitting the spectacular (but well documented) applications
to topology and strings.
Even so, it was impossible to render with mathematical rigor
the enormous amount of background material needed from the
physics.
I opted instead, in \S1, for the telling of one uniform but
very imprecise story of the physics.
My hope is that it will afford, to a mathematical reader
with no prior knowledge of the physics, at least a general
sense of the physical theory, of how the different parts
interact, and of what some of the basic concepts mean.
The remaining, mathematical, sections of these notes
discuss what is now known (and what is still missing) about
Seiberg-Witten type solutions of four-dimensional $N=2$
supersymmetric Yang-Mills theories in algebro-geometric
terms, with the emphasis on an interpretation via integrable
systems.

Section 1 starts with Maxwell's classical description of
electromagnetism.
We arrive at our central object of study, the duality in
super Yang-Mills theory, after discussing separately the
three directions in which Maxwell's theory needs to be
extended: nonabelian gauge groups, quantization, and
supersymmetry.
We are then ready for an overview of the Seiberg-Witten
solution, \cite{SW1}, emphasizing its algebro-geometric elements.
The section concludes with a discussion of the other
theories which could be exactly solvable by
algebro-geometric techniques.
These theories occupy us in the remainder of these notes.

Algebraically integrable systems are quickly reviewed in
\S2.
The upshot is that all the structures arising in the
Seiberg-Witten physics of \S1 are naturally realized here,
so algebraically integrable systems are the natural context
for describing the physics (cf. \cite{DW}, \cite{F}).
We also clarify the role of the Seiberg-Witten differential,
which occurs naturally in any integrable system of
{\it spectral} type.
One important feature of this differential, the linear
dependence of its residues on the mass parameters, turns out
to be part of a very general phenomenon.
We state and outline a proof of an apparently new result,
which is a complex analogue of the theorem of
Duistermaat-Heckman \cite{DH} on variation of the symplectic
structure in a family of symplectic quotients.
The complex case is somewhat subtler than the compact
version due to the presence of non-semisimple (e.g.
nilpotent) orbits.

Section 3 is mostly an attempt to rewrite the solution
\cite{DW} of the $\SU(n)$ theory with adjoint matter, with an
emphasis on the math rather than the physics.
We review the work of Markman \cite{Mn} and others on the
integrable systems of meromorphic Higgs bundles; among
these systems we pick the particular subsystem we need.
The system (i.e. its spectral curves) is described
explicitly, then compared to the elliptic soliton
description \cite{T, TV}.
A major part of \cite{DW} conisted of ``tests'' showing that
this system accurately reproduces physical properties of the
theory; we describe some of these ``tests'' here.

In the final section we discuss the other major class of
integrable systems which is known to solve the pure $N=2 $  
$\SYM$ theories: the periodic Toda lattices \cite{MW1}.
In particular, we compare the various forms assumed by the
spectral curves.
We then move to a common generalization of both systems,
the elliptic Calogero-Moser systems, which were
proposed in
\cite{Mc} as a solution of the theory with adjoint matter and
arbitrary gauge group.
This seems quite promising, but it still lacks a satisfactory
algebro-geometric description.
We also mention briefly some other theories where at least
the family of spectral curves is known, and conclude by
mentioning some recent ideas in string theory which exhibit
the $\QFT$ integrable systems discussed here as limits of
stringy systems.

It is a pleasure to acknowledge pleasant and useful
conversations, about some of the topics covered here, with
D. Freed, L. Jeffrey, A. Kirillov, D. Lowe, N. Nekrasov, R.
Plesser, N. Seiberg, S. Singer, R. Sjammar, J. Weitsman, and
E. Witten, as well as helpful comments on the manuscript by 
D. Freed, A. Grassi, A. Ksir, N. Nekrasov, M. Verbitsky 
and E. Witten.

\head
\S1. Some Physics Background
\endhead

\medskip\noindent
{\bf Electromagnetism}

The electric field $E$ and the magnetic field $B$ in
$\dbR^3$ are governed by Maxwell's  equations, of T-shirt
fame.
Think of $E$ and $B$, respectively, as a $1$-form and a
$2$-form on $\dbR^3$.
Both are time-dependent, so we combine them into the
electromagnetic field $F:= E\w dt+B$, a $2$-form on
$\dbR^4$.
Maxwell's equations in vacuo are:
$$
dF=0\qquad\qquad d*F=0.
$$
Here $*$ is the Hodge star operator on forms on $\dbR^4$,
determined by the Minkowski metric $g=dx^2+dy^2+dz^2-c^2
dt^2$.
(The speed of light $c$ is usually set to $1$ by adjusting
units.)

The same equations make sense for an unknown $2$-form $F$ on
any (pseudo) Riemannian manifold $(X,g)$.
They are clearly invariant under the group of global
isometries of $(X,g)$. 
When $\dim\,X=4$, the equations are also invariant under the
duality $F\mapsto *F$.
When $X=\dbR^4$ with the Minkowski metric, this duality
interchanges the electric and magnetic fields.
The group of isometries in this case is the Poincare group,
an extension of the Lorentz group $\SO(3,1)$ of rotations by
the group $(\dbR^4,+)$ of translations.

Geometrically, we fix a $U(1)$-bundle $L$ on $X$ and
consider a connection $A$ on $L$ as a new variable.
The first equation $dF=0$ is interpreted as saying that
(locally) $F=dA$ is the curvature.
The remaining equation then is $d*dA=0$.
The group of isometries of the triple $(X,g,L)$ 
acts on the connections and
preserves the equation.
This is an extension of the group of those isometries of
$(X,g)$ which preserve $L$ by the {\it gauge group} of
$C^\infty$ maps from $X$ to $U(1)$.
The duality is still there, but is now hidden.

Now impose an external field $j$.
(In $\dbR^4$, $j$ is a $3$-form whose $dx\w dy\w dz$
component is the charge density, while the other $3$
components give the current.)
The equations become
$$
dF=0,\qquad\qquad d*F=j,
$$
and duality is lost.
Geometrically, once we think of $F$ as curvature,  the first
equation becomes a tautology (the Bianchi identity), while
the second is an extra condition which may or may not hold.
Physically, we explain this asymmetry as due to the absence
of magnetic charges (``monopoles'').

Dirac proposed that some magnetic charges can be introduced
by changing the topology, or allowing singularities.
For example, a monopole at the origin of $\dbR^3$ is
represented by a $U(1)$-bundle on
$(\dbR_{x,y,z}^3\smallsetminus 0)\times\dbR_t$.
Such bundles are characterized topologically by their Chern
class $c_1\in\dbZ$, and this $c_1$ is interpreted as the 
magnetic charge at the origin. Indeed, the ``magnetic charge'' 
at the origin is just the total 
magnetic flux across a small sphere $S^2$ around the origin in 
$\dbR^3$; this is the integral over $S^2$ of $B$, which equals 
the integral of $F$, which is $c_1$. 
We can extend $F$  to $\dbR^4$ as a distribution, but now 
$dF$ is the
$\delta$-function $c_1.\delta_{x,y,z}$ instead of $dF=0$.
This is not yet sufficient to restore the symmetry to
Maxwell's equations: the electric charge is arbitrary, but
the magnetic charge is allowed only in discrete ``quanta''.
This does raise the possibility, though, that complete
symmetry could be restored in a quantum theory, where, 
as in the real world,
electric charge also comes in discrete units.
Indeed, Dirac discovered that in quantum mechanics, where 
an electrically charged particle is represented by a wave 
function satisfying Schrodinger's equation, the electric charge 
$e$ of one particle and the magnetic charge $g$ of another 
(in appropriate units) satisfy $e g \in \dbZ$ . This has the 
curious interpretation that all electric 
charges are automatically quantized as soon as one magnetic 
monopole exists somewhere in the universe.

In \cite{MO}, Montonen and Olive conjectured that in certain
supersymmetric Yang-Mills quantum field theories the
electric-magnetic duality is indeed restored.
This involves an extension of the basic picture in three
separate directions: replacing the ``abelian'' Maxwell
equations by the non-abelian Yang-Mills; quantizing the
classical theory; and adding supersymmetry to the Poincar\'e
invariance.

\medskip\noindent
{\bf Yang-Mills theory}

Yang-Mills theory is quite familiar to mathematicians.
It involves replacing the structure group $U(1)$ by an
arbitrary reductive group $G$.
We thus have a connection $A$ on some $G$-bundle $V$ over
our fourfold $X$, the curvature is now the $\ad(V)$-valued
two-form $dA+\frac{1}{2} [A,A]$, and the equations are as 
before,$dF=0$, $d*F=0$,
except that the covariant derivative $d=d_A$ on 
$\ad(V)$ now replaces the
exterior derivative $d$ on the trivial bundle $\ad(L)$.

These equations can be, and usually are, considered as
variational equations for the Lagrangian
$$
L(A):=\int_X \Tr(F\w *F),
$$
where the trace $\Tr$ can be any non-degenerate invariant 
pairing on the Lie algebra $\bfg$ of $G$.
This Lagrangian is considered as a functional on an
appropriate space of connections $A$ (if $X$ is not compact,
we must impose some decay condition at infinity for the
integral to converge).
The equation $dF=0$ is again the Bianchi identity, while
$d*F=0$ is the equation of the critical locus of $L$.
Either from the Lagrangian description or directly from the
equations, it is clear that the symmetry group includes
global isometries of $X$ which lift to $V$ (``the Poincare
group'') as well as local diffeomorphisms of the bundle
(``the gauge group'') and, in dimension $4$, the duality
$F\mapsto *F$.
In the original model of Yang and Mills, the group $G$ was
$\SU(2)$, while in QCD it is $\SU(3)$. 
For any ``realistic'' model of all of particle physics, $G$ would
need to contain at least $U(1)\times\SU(2)\times\SU(3)$.

\medskip\noindent
{\bf Quantization}

Quantization is an art form which, when applied to classical
physical theories, yields predictions of subatomic behavior
which are in spectacular agreement with experiments.
The quantization of the classical mechanics of particles
leads to quantum mechanics, which can be described, in a
Hamiltonian approach, via operators and rays in Hilbert
space, or in a Lagrangian formulation, via path integrals.
Quantization of continuous fields leads to QFT: quantum field 
theory. It has a Hamiltonian formulation, where the Hilbert space
is replaced by the ``larger'' Fock space, as well as a 
Lagrangian formulation, based on averaging the
action ($=\exp$ of the Lagrangian) over ``all possible
histories'' via path integrals.
These are generally ill-defined integrals over
infinite-dimensional spaces.
Nevertheless, there exist powerful tools for expanding them
perturbatively, as power series in some ``small'' parameters
such as Planck's constant, around a point corresponding to a
``free'' limit where the integral becomes
Gaussian and can be assigned a value.
(Actually, in modern QFT the ``small parameters'' used are 
often taken to be dimensionless coupling constants rather 
than Planck's constant.)
The coefficients of the expansion are finite dimensional
integrals encoded combinatorially as Feynman diagrams.
These power series typically diverge and need to be
``regularized'' and 
``renormalized''. These  processes are somewhat akin to the 
way we make the
Weierstrass $\grP$-function converge, by subtracting a
correction term from each summand in an infinite series:
a first attempt to write down a doubly periodic meromorphic 
function in terms of an infinite series of functions leads, 
unfortunately, to an everywhere divergent series; nevertheless 
we can obtain a series which does converge, and indeed to a 
doubly periodic meromorphic function, by subtracting an 
appropriate constant from each term. A physicist might say 
that we ``cancelled the infinity'' of the original series by 
subtracting ``an infinite constant''. In any case,
it is some of these renormalized values which are confirmed to
amazing accuracy by experiment, providing physicists with
unshakable confidence in the validity of the technique.

For a physicist, a ``theory'', either classical or quantum, is
usually
specified by its Lagrangian $L$, a functional given by
integration over the underlying space $X$ of a Lagrangian
density $\scrL$, which is a local expression involving the
values and derivatives of a given collection of fields.
The classical theory considers extrema of $L$, while the
quantum theory averages the action over all possible
fields.
The largest weight is thus still assigned to the extrema of 
$L$, but all ``quantum fluctuations" around the extrema are
now included. The Lagrangian density typically consists of a 
quadratic part, plus higher order terms which in a 
perturbative approach are interpreted 
as small perturbations. In a theory consisting of scalar fields
$\varphi_i$, the quadratic term might be 
$\sum _i (\vert d\varphi_i\vert^2+m_i^2 \varphi_i^2)$, where $m_i$
are the masses.
A perturbation might include terms like 
$g_I \varphi ^I$, 
where 
$I=(i_1, \dots, i_k)$
is a multiindex, and $g_I$ is a (supposedly small) 
{\it coupling constant} giving the strength of an interaction 
involving the $k$ fields $\varphi _i, i \in I$. 
(Other types of theories might involve also some fields which are 
vectors or spinors with respect to the Poincare group.) 
A typical quantity to be computed is the n-point function
$$ \langle x_1,\dots,x_n \rangle := 
\int \varphi (x_1) \dotsc \varphi (x_n) 
\exp{iL(\varphi)} D\varphi, $$
where integration is over the infinite-dimensional space of all 
$\varphi$'s, with respect to a (non-existent) measure $D\varphi$, 
and 
$L(\varphi) := \int _X \scrL (\varphi)$.

A formal approach to this integral is to interpret it as a power 
series in the coupling constants. Terms of the expansion are described 
by {\it Feynman diagrams}. The types of edges in a diagram correspond 
to the types of fields in the theory. The vertex types correspond to 
the multiindices $I$ with non-zero $g_I$. A $\varphi ^4$ term thus 
corresponds to a vertex where four $\varphi$ edges meet (this may 
represent two incoming particles which collide and scatter away), 
while a $\varphi \psi \bar{\psi}$ vertex can represent an electron 
$\psi$ annihilating a positron  $ \bar{\psi}$ to produce a photon 
$\varphi$. The {\it Feynman rules} assign to each diagram a $4\ell$ 
dimensional integral, where $\ell$ is the number of loops, i.e. first
Betti number, of the diagram, and 4 is the dimension of space-time: 
the integration is over all ways of assigning a momentum (a 4-vector) to 
each oriented edge, subject to Kirchhoff's law of null-sum at each 
vertex. The integrand for each diagram is the product over all edges 
of a basic function ({\it the propagator}) of the momentum assigned 
to that edge, times a factor involving the coupling constant $g_I$ 
at each vertex of type $I$, and an overall numerical factor accounting 
for graph automorphisms. The n-point function, or rather its Fourier 
tranform, is then expanded as the sum of all possible diagram 
integrals with $n$ external legs. 

One problem is that some of these Feynman integrals may be divergent. 
Of course, even if they converge individually, their sum may not. 
So the individual integrals need to be {\it regularized}, and the 
entire theory needs to be {\it renormalized}. A theory is said to be 
{\it superrenormalizable} if all divergences in it disappear after a 
finite number of its divergent diagrams are made finite by adjusting 
some of the parameters. In a {\it renormalizable} theory there are 
infinitely many divergent diagrams, but they can all be renormalized 
by iteratively adjusting the values of just a finite collection of
parameters. There is an easily verifiable necessary condition 
for renormalizability: each 
field $\varphi$ carries a certain dimensionality $[\varphi]$. (This 
is measured in powers of the mass; time and distance are converted 
to mass-inverse by setting the speed of light $c$ and Planck's 
constant $h$ to 1.) Now $L$, which occurs as an exponent, must 
have no dimensionality (it is a pure number, independent of units), 
so each term in $\scrL$ must have dimensionality $4$. This 
determines the dimensionality of each $g_I$. The obvious necessary 
condition for renormalizability  
is that $[g_I] \geq 0$  for all $I$: otherwise, 
the insertion of a vertex of type $I$ in an already divergent diagram 
will give a new diagram with divergence which is worse. 
Similarly, a necessary condition for superrenormalizability is that 
$[g_I] > 0$ for all $I$, since otherwise insertions into one 
divergent diagram would yield others with equally bad divergence.
These conditions serve to point out the very few 
candidates for renormalizable theories. Actually proving 
renormalizability in each case tends to be much harder.

For regularization, we introduce an auxiliary variable $\Lambda$ 
so that the integral can be interpreted as (the ``limit'' of) a 
function of $\Lambda$ going to infinity with $\Lambda$. For example, 
$\Lambda$ can be a momentum cutoff, meaning that the integral is 
carried out only over the ball of momenta $p$ satisfying
$\vert p \vert \leq  \Lambda$. (Or better, use some approximation of 
unity by smooth compactly supported cutoff functions.) 
Another popular variation is 
dimensional regularization: one writes down the expression for a 
given integral in a $d$-dimensional space, and observes that the 
answer is an analytic function in $d$ for sufficiently small $d$, 
typically acquiring a singularity at the relevant dimension $d=4$. 
The analogue of  $\Lambda$ in this version is $\exp{\frac{1}{4-d}}$.

The modern approach to renormalization is based on the {\it 
renormalization group flow}. Very crudely, one might illustrate 
this process as follows. The 
integral for each diagram $D$ translates into a function 
$f_D(m,g,p, \Lambda)$ of the masses, coupling constants, external 
momenta and  cutoff (or whatever else we used for $\Lambda$). If 
the number of divergent diagrams is finite (and $\leq$ the number
of parameters in $L$) then the simultaneous level sets of the 
corresponding $f_D$ will be (or will contain) a family of curves 
on which $\Lambda$ is unbounded. So we may hope to find a flow 
parametrized by $\Lambda$ along these curves:
$$ m=m(\Lambda), g=g(\Lambda), p=p(\Lambda).$$
In a general renormalizable theory there may be infinitely many 
$f_D$, but each may be modified (by a quantity which is bounded 
as a function of $\Lambda$) so there is still a flow tangent to 
all of them. Instead of attempting to fix the parameters at their 
``bare" values and taking $\Lambda$ to infinity, renormalization 
is accomplished by flowing along this renormalization group flow. 
The effective, physically measurable value of the parameters thus 
varies with the scale $\Lambda$. Similar phenomena are rather 
common in statistical physics. For example, the electric charge 
of a particle in a dielectric appears to be reduced, or 
{\it screened}, with distance, as a result of the alignment of 
the dielectric's charges around it. In QFT this occurs even in 
vacuum, because of the quantum fluctuations. 

A famous series of experiments at 
the Stanford Linear Accelerator showed that, when bombarded 
by very high energy electrons, a proton behaves as if it is made 
of three separate subparticles (quarks) which move freely within 
it. Yet at low energies (or equivalently, as viewed on a larger 
distance scale), the quarks are confined within the proton. A 
major success of QFT was its interpretation of these results in 
terms of the varying coupling constant $g$ for the strong force. 
The relevant theory (QCD: quantum chromodynamics) 
is {\it asymptotically free}, which 
means that $g \rightarrow 0$ as $\Lambda \rightarrow \infty$, 
and  $g \rightarrow \infty$ as $\Lambda  \rightarrow 0$. Thus 
at high energy scale $\Lambda$, the quarks behave as if $g=0$, 
i.e. as in the free theory, while at large distance 
$\Lambda \rightarrow 0$ the force becomes stronger and is 
presumably sufficient to confine the quarks. (This last point 
is not rigorously proved; Seiberg and Witten obtained a model 
for confinement in the supersymmetric analogue.)

An often computable quantity which describes the behavior of 
a coupling constant is its {\it beta function}, 
$\beta := \Lambda {\partial g/\partial \Lambda}$. An asymptotically 
free theory corresponds to $\beta <0$ for small values of $g$. 
In the boundary case 
where $\beta$ is identically  $0$, the theory is scale-invariant.

An additional complication is involved in the quantization
of theories which are (that is, their Lagrangian is)  
gauge-invariant. The problem with gauge-invariant theories 
stems from the degeneracy of their Lagrangian.
In classical gauge theory, the classical solutions (extrema of
the action) are not isolated, but come in entire gauge-group
orbits.
The equivalence of the Lagrangian and Hamiltonian
descriptions breaks down: this equivalence is based on a
Legendre transform, which makes sense only near isolated
points.
In the quantized theory, the path integral needs to be taken
not over all paths, but over gauge-equivalence classes.
This is seen already for the free (unperturbed) theory: the
value of a Gaussian integral is given by the inverse
determinant of the operator in the exponent; degeneracy (or
gauge invariance) implies that this determinant vanishes, and
needs to be replaced by a determinant on an appropriate
transversal slice or quotient. 
For quantized Yang-Mills, or gauge, theories, a good 
reference is \cite{FS}. The reward for handling these 
additional complications is a collection of theories which are 
renormalizable and sometimes asymptotically free: the beta-function 
of the basic theory (corresponding to Maxwell's in vacuo) is 
negative. The coupling to matter increases beta, but it remains 
negative when the number of quarks is small enough. 

\medskip\noindent
{\bf Supersymmetry}

Supermathematics is the habit of adding the prefix ``super''
to ordinary, commutative objects, to denote their
sign-commutative generalization.
Thus, a super space is a locally ringed topological space
whose structure sheaf $\scrO=\scrO_0\oplus \scrO_1$ is
$\dbZ_2$ graded, $\scrO_0$ is central, and multiplication
$\scrO_1\times\scrO_1\to\scrO_0$ is skew symmetric.
Affine super space $\dbR^{m,n}$ is the topologcal space
$\dbR^m$ with the sheaf $C^\infty(\dbR^m)\otimes
\Lambda^\cdot(\dbR^n)$ and the $\dbZ_2$ grading coming from the
$\dbZ$-grading of the exterior algebra.
An $(m,n)$-dimensional super manifold is a super space which
is locally isomorphic to $\dbR^{m,n}$ (The isomorphisms are,
of course, isomorphisms of $\dbZ_2$-graded algebras, and
need not preserve any $\dbZ$-grading.
A typical automorphism of $\dbR^{1,2}$, with even coordinate
$x$ and odd coordinates   $\theta_1$, $\theta_2$, may send
$x$ to $x+\theta_1\theta_2$.
So the identification of $\scrO$ as $C^\infty$ tensor
exterior algebras is not preserved under coordinate change.)
All of calculus extends to super manifolds.
The only non-obvious point is that the transformation law
for differential forms involves the determinant of the
transformation in the even variables, but the inverse
determinant of the transformation in the odd variables, and
a hybrid -- the {\it Berezinian} -- for a general transformation
mixing the parities.
The resulting operation of {\it integration} looks, in the odd
directions, more like usual {\it differentiation}, i.e., in local
coordinates $x_1,\dotsc,x_m,\theta_1,\dotsc,\theta_n$ it
reads off the coefficient of the top odd part, $\Pi_{i=1}^n
d\theta_i$.

The infinitesimal symmetries of a super manifold are
described by a super Lie algebra: a $\dbZ_2$-graded algebra
$\bfg=\bfg_0\oplus\bfg_1$ with a $[\,\,,\,\,]$ operation
which is sign-antisymmetric and satisfies the only
reasonable version of the Jacobi identity:
$$
\sum_{i=1}^3(-1)^{x_{i-1}x_{i+1}}[X_{i-1},[X_i,X_{i+1}
]]=0,
$$
where the $X_i$ ($i\in\dbZ_3$) are in $\bfg_{x_i}$.
In particular $\bfg_0$ is an ordinary Lie algebra, and
$\bfg_1$ is a $\bfg_0$-module with a $\bfg_0$-valued
symmetric bilinear form $[\,\,,\,\,]$ with respect to which
$\bfg_0$ acts as derivations:
$$
[X,[\theta_1,\theta_2]]=[[X,\theta_1],\theta_2]+
[\theta_1,[X,\theta_2]],
$$
for $X\in\bfg_0$, $\theta_i\in\bfg_1$.
The basic example is the Poincare super algebra
$\bfp=\bfp_0+\bfp_1$, whose even part $\bfp_0$ is the usual
Poincare algebra, and whose odd part $\bfp_1$ is the spinor
module $S$ of the Lorentz group $\SO(3,1)$.
(The Poincare group $P_0$ acts, via its quotient $\SO(3,1)$,
on $\bfp_1$; this induces the infinitesimal action of
$\bfp_0$.)
The pairing $\gamma\colon\,\bfp_1\times\bfp_1\to\bfp_0$ is
Clifford multiplication, with image the translation subgroup
of $\bfp_0$.
A physicist would write everything in coordinates, so
$\gamma$ becomes the collection of Pauli or Dirac matrices
(depending on whether real or complex forms are used).
The $N$-extended super Poincar\'e algebra has the same
$\bfp_0$, but $\bfp_1$ is replaced by $S^{\oplus N}$, the
sum of $N$ orthogonal copies of the spinors.

A supersymmetric (``$\SUSY$'') space is the super analogue
of a homogeneous space for the super Poincar\'e algebra or its
$N$-extended version.
An affine example is $(V,\scrO)$, where $V$ is a vector
space with non-degenerate quadratic form (of signature
$(3,1)$, in our case), and $\scrO$ is given in terms of the
spinor module $S=S(V)$ as
$C^\infty(V)\otimes\Lambda^\cdot(S^{\oplus N})$.
For more details on $\SUSY$ spaces, we refer to \cite{Be},
where they are defined to be super manifolds modelled on the
above affine example.

An $N=1$ $\susy$ theory is a theory (i.e. a Lagrangian
involving certain fields) with an infinitesimal action of
the super Poincar\'e algebra.
Likewise, a theory ``with $N$ supersymmetries'' has an
action of the $N$-extended Poincar\'e super algebra.
These theories can often, though not always, be described in
terms of a super Poincar\'e invariant Lagrangian involving
superfields on a $\susy$ space.
The ordinary Lagrangian, on ordinary space, is recovered by
Berezin integration along all the odd directions.
The classic reference on supersymmetry is \cite{WB}.
A possibly more friendly introduction is \cite{Fr}, and the
most geometrical version is \cite{Be}.
The first two references include a study of the
supersymmetric version of Yang-Mills theory.

In physics, there are two types of particles: bosons, 
which combine with each other freely, and
fermions, whose combinations are restricted by 
Pauli's exclusion principle. 
Mathematically, the distinction is roughly this:
consider a collection of $n$ particles with distinct
physical properties.
If $H_i$ is the quantum Hilbert space of states of the
i$^{\text{th}}$ particle, then the space of states of the
entire collection is $\otimes_{i=1}^n H_i$.
Now if instead we have $n$ identical particles, there are
fewer distinct states of the ensemble: for bosons we get 
$\Sym^n H$, while for fermions, $\Lambda^n H$.
In four dimensional theories, bosonic particles ordinarily
live in representations of the Poincare group with integer spin,
while fermions have half-integer (that is, non-integer) spin.
(This is the ``spin-statistics'' theorem.)
Most ``matter'' particles (electrons, protons, quarks) are
fermions; photons are bosons.

At a somewhat superficial level, the motivation for
introducing supersymmetry is to be able to treat the two
types of particles uniformly.
A complementary and deeper reason is based on the no-go
theorem of Coleman and Mandula
(cf. \cite{Fr}, \cite{WB}).
This says that the most general algebra of symmetries of a
class of quantum theories satisfying some rather reasonable
non-degeneracy assumptions (one of these is non-freeness,
another is the existence of massive
particles) is of the form $\bfp_0\oplus \bfg$, where
$\bfp_0$ is the Poincar\'e algebra of space-time isometries,
and $\bfg$ is the gauge algebra.
The point is that this is a direct sum (of algebras), so
there can be no non-trivial mixing of the two symmetries.
Since symmetries are equivalent to conserved quantities,
this implies for instance that any conserved quantity
(``current'') other than energy-momentum and angular
momentum (which come from $\bfp_0$) must transform as a
scalar under $\bfp_0$.
This is unsatisfactory, since some free theories (e.g. of
one real scalar field plus one real vector field, 
\cite{Fr}; such
theories are not covered by the Coleman-Mandula theorem) do
admit conserved vectors and tensors, which should survive
under small perturbation.

The tasks of mixing the two particle types and of avoiding
the Coleman-Mandula restriction can fortunately be
accomplished simultaneously through the introduction of
super Lie algebras: even elements preserve particle types,
odd elements exchange them.
The even part of the algebra obeys Coleman-Mandula, but
there is room for odd symmetries and the corresponding odd,
or spinorial, conserved quantities.
A theorem of Haag, Sohnius and Lopuszanski classifies the
possible super Lie algebras of symmetries of physically
interesting theories: there are the $N$-extended super
Poincare algebras, as well as some central extensions and
variants including an (even) gauge algebra $\bfg$ acting on
the odd sector.

\medskip\noindent
{\bf N=2 Super Yang-Mills}

The setting for Montonen-Olive duality \cite{MO}, for the work
of Seiberg-Witten \cite{SW1, SW2}, and for the development
discussed in later sections, is $N=2$ supersymmetric
Yang-Mills theory ($\SYM$) in four-dimensional space.
The Lagrangian for this theory is usually written in terms
of super fields on an affine SUSY space \cite{WB}, making
the super Poincare invariance evident.
When written out explicitly in terms of its component
fields, the Lagrangian is quite complicated.
It starts with a purely bosonic term analogous to the usual
Yang-Mills Lagrangian $\int\Tr(F\w*F)$.
Additional terms, involving additional fields, are required
for the $N=2$ supersymmetry.
The fields involved are grouped into $N=1$ multiplets
(representations of the $N=1$ super Poincare algebra), which
in turn combine into $N=2$ multiplets.

Such a theory depends, of course, on the choice of a compact
gauge group $G$.
For a given $G$ there is the pure $N=2$ theory, analogous to
Maxwell's equations in vacuo, containing only those fields
and terms in the Lagrangian required for supersymmetry.
There are also various theories in which some combination of
additional fermionic (``matter'') particles is thrown in.
We will return to these shortly.

Each $\SYM$ theory has a moduli space $B$ of vacuum states
($=$ eigenstates of the Hamiltonian or energy operator 
corresponding to the lowest eigenvalue).
This is of course a consequence of the degeneracy, or of the
gauge invariance, of the Lagrangian: a non-degenerate Lagrangian 
should correspond to a unique vacuum, and an ordinary degeneracy 
is obtained when two eigenvalues happen to coincide, so there are 
two (or more) independent vacua. The situation in gauge theory 
is that there is a continuum $B$ of independent vacua. The 
complex dimension of $B$ is the rank $r$ of $G$, and we
can in fact describe a coordinate system $(u_i)$,
$i=1,\dotsc,r$, on it.
First, the classical analogue: the classical potential is
$$
V(\varphi)=\frac{1}{g^2}\Tr[\varphi,\varphi^\dagger]^2,
$$
where $g$ is a coupling constant (a parameter in the
Lagrangian), and $\varphi$ is a Higgs field -- a complex
field transforming as a scalar under Poincare, and in the
adjoint representation of $G$.
Since $[\varphi,\varphi^\dagger]$ is self-adjoint, the minima
occur precisely when $\varphi$ commutes with $\varphi^\dagger$.
This locus is $G$-invariant: it is the union of the $G$-orbits 
of elements in the Cartan subalgebra $\bft$ of the complexified 
Lie algebra  $\bfg_{\dbC}$. 
The quotient by $G$ is thus $B=\bft/W$,
where $W$ is the Weyl group. 
For example, when $G=\SU(n)$ we have $\bft\approx
\dbC^{n-1}$, and the quotient is again $\dbC^{n-1}$, with
coordinates $u_2=\Tr\,\varphi^2,\dotsc,u_n=\Tr\,\varphi^n$.
(The reader will no doubt recognize this as a (real) 
symplectic reduction: the action of $G$ on $\bfg_{\dbC}$ is 
Hamiltonian with respect to the natural translation-invariant 
symplectic structure induced on  the real vector space 
underlying $\bfg_{\dbC}$ by the 
imaginary part of the Killing form. The moment map sends 
$\varphi$ to $\Tr[\varphi,\varphi^\dagger]$, and the 
symplectic quotient is  $\bfg_{\dbC} // G = B = \bft /W$.
We will run into the complex analogue of symplectic 
quotients in the discussion of linearity in the next section.)

Now the constraints imposed by $N=2$ supersymmetry (i.e. the
difficulty of writing down a supersymmetric Lagrangian)
imply that the quantum moduli space remains the same $B$, as
an abstract algebraic variety, so it still has the
coordinates $u_i$.
But we will see that just about all other features of $B$
change in going from the classical to the quantum: 
for example, this applies to the 
K\"ahler metric on it, the locus of singularities (of the
metric), the flat structure, and the relation of the $u$
coordinates to the natural fields in the theory.

A central quantity in any $\SYM$ theory is the electric
charge $a$, which is the vector of eigenvalues of the Higgs
field $\varphi$.
It lives in the Cartan, $\bft$.
For $\SU(2)$ this is a complex number, as is $u$.
We can take $\varphi$ to have eigenvalues $\pm\,\frac a2$,
and then $u=\Tr\,\varphi^2=\frac{a^2}{2}$, in the classical
theory.
In the quantum theory, this relationship still holds
approximately as $u\to\infty$, by the asymptotic freedom of
the theory.
But for finite $u$ we will see that the relationship is very
different.

A crucial property of these theories is that their
low-energy (read: real-world observable) behavior can be
described in terms of only a finite number of parameters --
the coordinates $u_i$ on the quantum moduli space $B$.
This description is via a ``low energy effective
Lagrangian'' $\scrL_{\eff}$, obtained from the full
Lagrangian through a process of averaging over all heavy
degrees of freedom.
So all measurable quantities become (locally) functions on
$B$.
In particular, this applies to the components $a_i=a_i(u)$
of the electric charge.
Supersymmetry implies that these functions are holomorphic,
and we will use them as alternative local coordinates on
$B$.

While the existence of $\scrL_{\eff}$ is guaranteed by the
general theory, its actual computation can be extremely
difficult.
The supersymmetry allows $\scrL_{\eff}$ to be expressed
locally in terms of a single holomorphic function
$\scrF=\scrF(a)$, the prepotential.
(This is due to Seiberg, cf. \cite{SW1}, (2.7) and (2.11).
Note that we need $\scrF$ as a function of $a$, not of $u$.)
All other quantities in the theory can then be expressed in
terms of $\scrF$ and the coordinates $a$.
This includes the vector
$$
a^D=\frac{d\scrF}{da},
$$
whose physical significance will appear shortly; the matrix
$$
\tau(a)=\frac{d^2\scrF}{da^2}=\frac{da^D}{da},
$$
which is a complexified coupling constant (for $\SU(2)$, the
low-energy  effective values of the real coupling constant
$g$ (encountered in the potential $V(\varphi)$) and of a
phase angle $\theta$ are given by
$\tau=\frac{\theta}{2\pi}+\frac{4\pi i}{g^2}$.
In particular, $\scrF(a)$ must be such that $\im(\tau)>0$);
the K\"ahler metric
$$
ds^2=\im(\tau\,da\,d\abar)=\im(da^D\,d\abar)
$$
which was mentioned above; and its K\"ahler potential 
$$
K=K(a,\abar)=\im\left(\sum_i a_i^D\abar_i\right)
=\im\left(\sum_i \frac{d\scrF}{da_i}\,\abar_i\right).
$$
The effective Lagrangian itself is given as a sum of kinetic
and potential parts, each a Berezin integral over
an $N=1$ SUSY affine space 
of dimensions (4,4) and (4,2), respectively:
$$
\scrL_{\eff}=\frac{1}{4\pi}\int K(\Phi,\bar{\Phi})d^4\eta+
\frac{1}{8\pi}\im\int\tau(\Phi)\vert W\vert^2 d^2\eta.
$$
Here $\Phi$ and $W$ are the two $N=1$ multiplets which together
make up the $N=2$ multiplet containing the Higgs field
$\varphi$.
($\varphi$ and a fermion $\psi$ are the fields in the
``chiral multiplet'' $\Phi$. The ``vector multiplet'' $W$
contains another fermion $\lam$ and the gauge field $A$.
It can be thought of as the field strength, $N=1$ super
analogue of the curvature $F$ in ordinary $YM$.)
The integrands involve an analytic function ($\tau$)
or the real part of one ($K$); this function 
is to be applied formally to the super
fields $\Phi$, $W$ and their conjugates.
The integration over the four (or two) odd space directions 
$\eta$ is interpreted via Berezin's recipe as a derivative, 
or as the reading off of the relevant coefficients.
(The second term is actually written
as an integral over a (4,2)-dimensional ``chiral'' quotient of the full
(4,4)-dimensional SUSY space, in which two of the odd directions have
been reduced.

What is missing is an explicit expression for $\scrF$.
Traditionally, $\scrF$ was written as an infinite sum
involving its classical (``tree level'') value, the one-loop
correction which is logarithmic in $a$ (higher loop
contributions were known to vanish), and an infinite sum of
mysterious instanton corrections involving negative powers
of $a$.
The breakthrough came when Seiberg and Witten succeeded, in
\cite{SW1}, in computing $\scrF$ (and hence $\scrL_{\eff}$,
$ds^2$, etc.) exactly, for the pure $N=2$ theory with
$G=\SU(2)$.
Their solution is based on the global properties of the
theory.

\medskip\noindent
{\bf Duality}

The original conjecture of Montonen and Olive \cite{MO} was
that the low energy theory should be invariant under a
duality exchanging electric and magnetic charges, and also
exchanging the coupling constant $g$ with its inverse.
In other words: electric charge in the strongly coupled theory
(large $g$) behaves like magnetic charge at weak coupling,
and vice versa.
Seiberg and Witten found that a modified version of this
duality holds in pure $N=2$ $\SYM$ for $\SU(2)$.
They applied the super version of the Hodge $*$ operator to
the Lagrangian $\scrL_{\eff}$, added some variables
(Lagrange multipliers) and integrated out others, and ended
up with a new Lagrangian $\scrL_{\eff}^D$ which is
equivalent to the original, and has the same form, but in
terms of a new quantity $a^D:=\frac{d\scrF}{da}$:
$$
\scrL_{\eff}^D(a^D)=\scrL_{\eff}(a).
$$
The same transformation sends $a^D$ to
$-a=-\frac{d\scrF}{da^D}$.
(This must be so, since the first term in $\scrL_{\eff}$
changes sign when $a$, $a^D$ are simply exchanged.) 
The {\it physical meaning} of $a^D$ is therefore that of the dual,
or magnetic, charge.
As to the coupling constants, the defining relations
$$
\tau=\frac{da^D}{da},\quad
\tau^D=\frac{d(-a)}{da^D}
$$
together with the chain rule, imply that the actual duality
transformation is
$$
\tau^D=\tau^D(a^D)=-\tau(a)^{-1}.
$$
This specializes to the Montonen-Olive transformation
$$
g^D=g^{-1}
$$
when $\theta=0$, but not otherwise.

There is another, more elementary, transformation which
preserves the Lagrangian: the $\theta$ angle is defined only
modulo $2\pi\dbZ$, so $\tau$ is defined modulo $\dbZ$.
There is therefore a transformation which fixes $a$ and
sends $\tau\mapsto\tau+1$.
Since $\tau=\frac{da^D}{da}$, this requires that $a^D$ goes
to $a^D+a$.

The two transformations can be represented by the matrices
$$
S=\pmatrix\format \r &\quad \r\\
0 &1\\
-1 &0
\endpmatrix
\qquad\qquad
T=\pmatrix
1 &1\\
0 &1
\endpmatrix\,\,,
$$
acting linearly on the $2$-vector $\left(\smallmatrix
a^D\\ a\endsmallmatrix\right)$ and fractional-linearly on
$\tau$.
Since $S$ and $T$ generate $\SL(2,\dbZ)$,we see that the
theory has an $\SL(2,\dbZ)$ action.
(For other groups, $\SL(2,\dbZ)$ is replaced by a subgroup
of finite index in $\Sp(2r,\dbZ)$.)

\medskip\noindent
{\bf The Seiberg-Witten solution}

The foregoing suggests that what lives intrinsically over a
generic $u\in B$ is not the electric charge $a(u)$ but the
unimodular lattice $\dbZ a(u)+\dbZ a^D(u)$ of all charges
(electric, magnetic, and ``dyonic'', or mixed).
As $u$ varies, we get a $\dbZ^2$ local system $V$ over $B$
minus some singular locus, together with a homomorphism to
(the trivial bundle with fiber) $\dbC$.
Dually, this homomorphism is equivalent to a holomorphic
section $\left(\smallmatrix a^D\\ a\endsmallmatrix\right)$
of the vector bundle $V_{\dbC}$ (i.e. of the sheaf
$V\otimes\scrO_{B\smallsetminus\sing}$).
The plan is to identify the local system (variation of Hodge
structure, really) $V$, and then to fix the section
$\left(\smallmatrix a^D\\ a\endsmallmatrix\right)$, which in
turn will determine the prepotential $\scrF$ (up to
constant).

Seiberg and Witten identify $V$ by finding its monodromy.
Near $u=\infty$ we know that $u\sim\frac12\,a^2$, so the
monodromy around $\infty$ sends $a\mapsto -a$.
The known non-trivial logarithmic (one-loop) term in the
perturbative expansion of $\scrF$ then fixes the shift in
$a^D$: the monodromy is 
$$
M_\infty=
\pmatrix\format \r &\quad\r\\
-1 &2\\
0 &-1
\endpmatrix\,\,.
$$
In particular, the local system is non trivial, so $\tau$ is
a non-constant holomorphic function on the universal cover
of $B\smallsetminus\sing$, with values in $\dbH$.
The singular locus must therefore contain at least two
points (in addition to $\infty$). From a physical argument,
Seiberg and Witten deduce that the system is as simple
as possible: exactly three singularities, which can be
relocated to $u=\pm1$ and $\infty$ by an automorphism.
Duality instructs that the finite monodromies should be
conjugate to $-M_\infty$.
Seiberg and Witten 
label things so that $a^D$ vanishes at $u=1$, yielding
the monodromy
$$
M_1=\pmatrix\format \r &\quad\c\\
1 &0\\
-2 &1
\endpmatrix
$$
which then determines $M_{-1}$.

The monodromies generate the level-$2$ congruence subgroup
$\Gamma(2)\subset \SL(2,\dbZ)$, and the local system itself
can now be identified with the fiber cohomology of the
universal level-$2$ elliptic curve
$$
E_u\colon\,y^2=(x+1)(x-1)(x-u)
$$
over the $u$-line $\dbC\smallsetminus\{\pm1\}$.

The complexification $V_\dbC$ has a global trivialization in
terms of the holomorphic $1$-form $\lam_1=\frac{dx}{y}$ and
the residueless meromorphic form $\lam_2=\frac{xdx}{y}$.
We can also fix an integral homology basis, say
$\gamma=\Loop$ around the branch points $0,1$, and
$\gamma^D=\Loop$ around $1$, $u$.
To determine $a(u)$ and $a^D(u)$, we use the condition that
$\im(\tau)>0$, where $\tau(u):=\frac{da^D/du}{da/du}$.
One geometric solution is to take $\tau(u)$ as the period
$$
\tau_u:=
\frac{\oint_{\gamma^D}\lam_1}{\oint_{\gamma}\lam_1}
$$
of the elliptic curve $E_u$; rigidity implies that this is
in fact the only solution.
We then find that, up to a multiplicative constant, $a$ and
$a^D$ are the periods over $\gamma$, $\gamma^D$ of the
meromorphic $1$-form
$$
\lam=\frac{y\,dx}{x^2-1}=\frac{(x-u)dx}{y}=\lam_2-u\lam_1.
$$

\medskip\noindent
{\bf Adding matter}

In a quantum field theory, a matter particle is represented
by a fermion field $\psi$ which, in the case of a gauge
theory, needs to be in some representation of $G$.
It is added (or: coupled) to the theory by adding to the YM 
Lagrangian a term corresponding to the Lagrangian of the
free particle (this includes a parameter, the particle mass
$m$), and another term corresponding to its interaction with
the YM field.
In $\SYM$ the field needs to be in a representation $\rho$
of $G$, tensored with an $N=2$ multiplet, to preserve the
$\susy$.
(The relevant multiplet, which is built from the
fermion $\psi$, is called a {\it hypermultiplet}; we have
already encountered the other basic type, the $N=2$ {\it
vector multiplet}, which contained the scalar Higgs field
$\varphi$.)
If the number of particles added is not too large, the
theory remains asymptotically free $(\beta<0)$ or becomes
scale invariant $(\beta=0)$, and in either case retains many
properties of the pure $\SYM$ theory.
There still is an $r$-dimensional moduli space $B$.
At $u\in B$ there are $r$ electric charges $a_i(u)$ and $r$
dual, magnetic, charges, $a_i^D(u)$, each charge living in
the Cartan $\bft\approx \dbC^r$ as before.
Duality is modified somewhat: instead of $\Sp(2r,\dbZ)$
transformations acting on $\left(\smallmatrix a^D\\
a\endsmallmatrix\right)$, one gets affine symplectic
transformations of the form
$$
\pmatrix
a^D\\
a\endpmatrix\longmapsto R\pmatrix
a^D\\
a\endpmatrix +\sum_{i=1}^{N_f} m_i
\pmatrix
n_i^D\\
n_i\endpmatrix
$$
where $R\in\Sp(2r,\dbZ)$, the $m_i$ are the masses of the
$N_f$ particles added, and $n_i$, $n_i^D$ are integral
$r\times r$ matrices.
(``$N_f$'' stands for number of {\it flavors}.
The number of {\it colors} $N_c$ is the dimension of the
representation space of $\rho$.
The terminology comes from quantum chromodynamics (QCD), a
gauge theory with $G=\SU(3)$.)

It is reasonable to guess that this structure can be
implemented by a family of $r$-dimensional abelian varieties
$A_u$ over $u\in B\smallsetminus(\sing)$:
The $r\times 2r$ period matrix of $A_u$ in some integral
homology basis $\gamma^D$, $\gamma$ is 
$\left(\frac{da^D(u)}{du},\frac{da(u)}{du}\right)$.
(The polarization on the complex torus corresponding to the
lattice generated by columns of this matrix is given by
the Dirac quantization condition.)
On the total space of the family there is a meromorphic
$1$-form $\lam$, the {\it Seiberg-Witten differential}.
There are $N_f$ divisors $D_i$ along which $\lam$ acquires a
pole with constant residue $\frac{m_i}{2\pi i}$.
The charges $a^D$, $a$ can then be recovered as the periods
of $\lam$ over $\gamma^D$, $\gamma$.
These have precisely the ambiguity needed: monodromy of the
cycles $\gamma^D$, $\gamma$ acts as $\Sp(2r,\dbZ)$, while
moving such a cycle across the divisor $D_i$ changes the
period by an integer multiple of $m_i$.
The masses $m_i$ acquire a cohomological interpretation:
the $2$-form $\sigma=d\lam$ is holomorphic (in fact, it is a
holomorphic symplectic form) on the total space, and its de
Rham cohomology class can be expressed in terms of the
fundamental classes of the $D_i$:
$$
[\sigma]=\sum m_i[D_i].
$$
In particular, as we vary the masses, we obtain a family of
algebraic symplectic manifolds in which the symplectic class
varies {\it linearly} in the parameters.

In the case of scale-invariant theories ($\beta=0$), the
theory and therefore also its geometric model depend on one
parameter in addition to the masses.
This is $\tau_{\cl}$, the classical limit $(u\to\infty)$ of
the quantum coupling constant $\tau=\tau(u)$.
In asymptotically free theories $(\beta<0)$ we have
$$
\lim\limits_{u\to\infty}\tau(u)=\infty,
$$
and instead of $\tau$ the theory depends on a scale
parameter $\Lam$.
In the \cite{SW1} solution, this $\Lam$ is set to $1$ for
convenience.

In \cite{SW2} Seiberg and Witten show that the above
geometric model indeed leads to an exact solution of the
various $\SU(2)$ $\SYM$ theories with matter.
For each of these, the solution requires a family
parametrized by the vector $m=(m_i)$ of masses (and by the
additional parameter $\tau$ or $\Lam$, which we omit from
notation), of elliptic surfaces 
$$
X_m\to B\approx\dbC
$$
with meromorphic $1$-forms $\lam_m$ whose residues equal
$\frac{m}{2\pi i}$.
The solution satisfies numerous consistency checks (the way
one theory is known to flow to another in the limit where
some parameter, such as a mass, goes to $0$ or $\infty$, is
reflected in corresponding degenerations of one geometric
model to another) and has all the internal symmetry of each
of the physical theories.

The number of particles which can be added to a SYM theory
is determined as follows.
The Killing form on the reductive Lie algebra $\bfg$
determines an element $c_2$ in the center of the universal
enveloping algebra $U(\bfg)$.
By Schur's lemma, $c_2$ acts on the total space of each
irreducible representation $\rho$ as multiplication by
some scalar $c_2(\rho)$.
(This $c_2$ is closely related to the second chern class 
of a bundle: 
 Let $V$, $V'$ be two
vector bundles on a surface associated via representations
$\rho$, $\rho'$ to the same $G$-bundle.
Then
$$
\frac{c_2(\rho')}{c_2(\rho)}=\frac{c_2(V')}{c_2(V)},
$$
where $c_2$ on the right denotes Chern classes. This can 
be used to give another definition of $c_2$ of a representation.)
In any case, the $\beta$ function of the SYM theory with
particles in representations $\rho_i$, $i=1,\dotsc,N_f$, is
$$
\beta=c_2(\ad)-\sum_i c_2(\rho_i)
$$
where $\ad$ is the
adjoint representation.
For the classical algebras we can adjust the scaling of the
Killing form so that $c_2(F)=1$ where $F$ is the fundamental
representation.
With this normalization, $c_2(\ad)=2h_\bfg^{\upvee}$, where
$h_\bfg^{\upvee}$ is the dual Coxeter number.
For $SU(n)$ we have $h_\bfg^{\upvee}=n$, so the $\SU(2)$
theory, for example, remains asymptotically free
with$N_f=0,1,2$ or $3$ particles in the fundamental
representation $F$, and becomes scale invariant for $N_f=4$
fundamental particles or for one particle in the adjoint
representation.
These are precisely the cases studied and solved in
\cite{SW2}.

It is natural to ask whether the four dimensional, $\ssize N=2$ SYM
theory with arbitrary gauge groups $G$ and collection of
particles $\rho_i$ with
$\beta\Le 0$ can be similarly solved in
geometric terms via families of abelian varieties and
meromorphic one-forms.
This question will occupy us in the remainder of this work.

\head \S2.  Why integrable systems
\endhead

\medskip\noindent
{\bf Algebraically integrable systems}

An algebraically integrable system consists of a complex 
algebraic manifold $X$, an everywhere non-degenerate, closed
holomorphic $(2,0)$-form $\sigma$ on $X$, and a morphism 
$\pi : X\to B$ whose
general fiber is a (polarized) abelian subvariety 
$X_b: = \pi^{-1} (b)$ which is
Lagrangian with respect to $\sigma$.  In practice, the object 
which is naturally
given is often an affine open $X_0\subset X$, the fibers being 
affine subsets of
abelian varieties.  This is a straightforward algebraic analogue 
of the (real,
$C^\infty$) notion of integrable system in classical mechanics.  
As there, any
function $H$ (``Hamiltonian'') on the base $B$ determines a 1-form 
$d\pi^*H$ on
$X$, and when contracted with the symplectic form $\sigma$ it gives 
a vector
field $v_H$ on $X$, tangent to the fibers of $\pi$.  The pullbacks 
$\pi^*H_1,
\pi^*H_2$ of two functions on $B$ Poisson-commute, and in particular 
the
corresponding vector fields commute.

Locally on the base of an algebraically integrable system 
there is a natural holomorphic flat
structure, depending on a discrete choice.  The fiber homology 
$H_1(X_b, 
\dbZ)$ has a non-degenerate skew-symmetric pairing (induced 
by the 
polarization).  The choice involved is that of a Lagrangian 
subspace $V^*_{\dbZ}$
for this pairing.  This is often achieved by specifying a 
symplectic basis 
$\alpha_1,\dotsc, \alpha_g, \beta_1,\dotsc, \beta_g$ of 
$H_1(X_b,\dbZ)$, and
taking the subspace $V^*_{\dbZ}$ spanned by the $\alpha$ 
cycles.  This
choice needs to be made for one $b\in B$, and (by Gauss-Manin) 
is propagated 
to nearby fibers.  Integration over these $\alpha$-cycles gives 
a local
trivialization of the relative cotangent bundle $\pi_*\omega_X$, 
and the
symplectic form $\sigma$ converts this to a flat structure on 
(the tangent bundle
of) $B$, modelled on $V: = Hom(V^*_{\dbZ},\dbC)$.

Locally on the flat base $B$, we can ask what kind of data is 
required to specify
an algebraically integrable system inducing the given flat 
structure.  We think
of the family of Abelian varieties, plus choice of $\alpha$ 
and $\beta$ cycles,
as specified by its period map 
$$
p:B\to \dbH_g
$$
to the Siegel half space $\dbH_g$ of $g\times g$ symmetric 
complex matrices
with positive imaginary part.  The total space $\pi:X\to B$ 
is then retrieved
as pullback of the universal abelian variety 
$\scrX{_g}\to\dbH_g$.  The 
question then is to understand the condition on $p$ for 
there to be a symplectic
form on $X$ for which $\pi$ is a Lagrangian fibration.  
The fixing of the
flat structure is achieved by thinking of $\dbH_g$ as 
open subset of 
$Sym^2 V^*$, where $V$ is the tangent space to $B$ 
(at any $b\in B$).  The
derivative of the period map is
$$
dp_b\in Hom(V, Sym^2 V^*) = V^*\otimes Sym^2 V^*.
$$
A straightforward computation done in \cite{DM2} shows that 
the condition is that
$dp_b$ lives in the subspace Sym$^3 V^*$, for all $b\in B$.  
So an integrable
system with base $B$ is specified by a field of cubics on the 
tangent bundle of
$B$.

This field of cubics is subject, of course, to a strong 
integrability condition.
In flat local coordinates $z_i$ on $B$, let $(p_{ij}(z))$ 
be the period matrix
and $(c_{ijk}(z))$ the field of cubics.  Then equality of 
mixed partials
$$
\frac{\partial p_{ij}}{\partial z_k} = c_{ijk} =
\frac{\partial p_{ik}}{\partial z_j}
$$
implies the existence of functions $w_i$ (locally) 
on $B$ such that
$$
p_{ij} = \frac{\partial w_i}{\partial z_j}.
$$
These again have symmetric partials,
$$
\frac{\partial w_i}{\partial z_j} = p_{ij} =
\frac{\partial w_j}{\partial z_i},
$$
so there is a single function $\scrF$ on $B$, called 
the {\it prepotential},
satisfying:
$$
\align
w_i &= \frac{\partial\scrF}{\partial z_i}\\
p_{ij} &= 
\frac{\partial^2\scrF}{\partial z_i \partial z_j}\\
c_{ijk}
&= \frac{\partial^3\scrF}{\partial z_i\partial z_j\partial z_k}
\endalign
$$
The prepotential determines also a Kahler potential on $B$
$$
K: = Im(\sum_i w_i\barz_i) = 
Im(\sum_i\frac{\partial \scrF}{\partial z_i}\barz_i)
$$
and its associated Kahler form
$$
\frac i2\partial\overline\partial K = \frac i2
\sum_{i,j}Im(p_{ij}) d z_i\wedge d\barz_j.
$$
Equivalently, this Kahler (1,1)-form on $B$ is obtained by 
integration
over the $g$-dimensional fibers $X_b$ of the
$(g+1, g+1)$-form $\sigma\wedge\barsig\wedge t^{g-1}$ on $X$, 
where $t$ is any
$(1,1)$-form on $X$ whose restrictions to fibers give the 
polarization.

Conversely, given $B$ with coordinates $z_i$ and a holomorphic 
function $\scrF$,
the above formulas produce an algebraically integrable system 
over the open
subset of $B$ where 
$Im(\frac{\partial^2\scrF}{\partial z_i\partial z_j})>0:$
the family of abelian varieties is determined by its 
periods $p_{ij}$, and the
symplectic form is given by the cubic $c_{ijk}$.  The 
slightly weaker notion of
an analytically integrable system (just replace abelian 
varieties by polarized
complex tori) is recovered over the open subset where
$Im(\frac{\partial^2\scrF}{\partial z_i\partial z_j})$ 
is invertible, cf. 
\cite {DM2}.  The Kahler ``metric'' in this case is still 
non-degenerate, but
possibly indefinite.

We note that the $w_i$ play a role dual to that of the 
$z_i$: the $z_i$ 
give the flat structure on $B$ determined by integration 
over the $\alpha$ cycles,
while the $w_i$ give the flat structure determined by 
integration over the $\beta$
cycles.  We also note that, up to this point, the $z_i$ 
as well as the $w_i$ have
been determined only up to arbitrary additive constants.  
In other words, what we
have are the vector fields 
$\partial/\partial z_i, \partial/\partial w_i$, or the 
1-forms $dz_i, dw_i$.  We
will discuss below the additional choice used in the SW setup 
to fix these
additive constants.

Similar structures have arisen in various more specialized 
contexts, e.g. in the
WDVV equations of conformal field theory and on moduli spaces 
of Calabi-Yau
threefolds, as well as in the Seiberg-Witten setup described 
in \S 1.  They often
go under the title of ``special geometry''.  The point is that 
these specialized
contexts are not really needed:  the structure of special 
geometry on $B$
(consisting of a covering by open subsets on which an 
appropriate prepotential,
dual flat holomorphic coordinate systems, and a Kahler 
metric are defined),
arises naturally on the base of any algebraically 
integrable system.  The
converse is also true, locally, as we have just seen.  
Globally, though, there is
an integrability constraint:  the monodromy action 
(say on the coordinates $z_i,
w_i$) in an integrable system involves integral 
symplectic transformations plus
(complex) translations, while special geometry 
allows the symplectic
transformation to be real instead of integral.  
We refer the reader to \cite{F}
for a very clear differential-geometric discussion 
of special geometry and its
relationship to integrable systems.

\medskip\noindent
{\bf Seiberg-Witten Differentials}

In the Seiberg-Witten picture of super Yang-Mills, 
there are the charges 
$a_i$ and $a_i^D$which are described locally as holomorphic 
functions on
the quantum moduli space $B$.  There is less ambiguity in 
choosing these $a_i,
a_i^D$ than in the coordinates $z_i, w_i$ on the base of 
an algebraically
integrable system:  the group of their linear combinations 
$\Sigma n_i a_i +
n_i^Da_i^D$ (with integer coefficients $n_i$) should be 
uniquely determined in
the pure SYM case, and determined modulo the masses of the 
added particles in
general.  A natural way to get functions with precisely 
such ambiguity on the
base of an integrable system $\pi: X\to B$ is to choose a 
meromorphic
differential 1-form $\lambda$ on $X$ and let $a_i, a_i^D$ 
be its periods over a
set of 1-cycles representing $\alpha_i,\beta_i$ respectively 
and avoiding the
poles of $\lambda$.  We want $a_i, a_i^D$ to be a possible 
choice of the
coordinates $z_i, w_i$, which are defined up to translation 
(once we have fixed
the $\alpha_i, \beta_i$).  In other words, we want
$$
da_i = dz_i, \quad da_i^D = dw_i.
$$
Since $dz_i, dw_i$ were defined as the contraction of 
$\alpha_i, \beta_i$
respectively with the symplectic form $\sigma$, the condition 
becomes simply:
$$
d\lambda = \sigma.
$$
The behavior of the singularities of $\lambda$ is constrained 
by this condition.
When restricted to a general fiber $X_b$, $\lambda$ will have 
poles along the
union of some irreducible divisors $D_{b,j}$.  One constraint
is that the
residue Res$_j$ of $\lambda$ along $D_{b,j}$ is (locally in $B$) 
independent of
$b$.  A {\it Seiberg-Witten differential} is a meromorphic 1-form 
$\lambda$
satisfying $d\lambda = \sigma$ and Res$_j = m_j$, the mass of the 
$j$-th
particle.  On an algebraically integrable system with a specified 
Seiberg-Witten
differential, the local coordinates $z_i, w_i$ can be specified as
the charges $a_i, a_i^D$, with just the right ambiguity.

One way to describe an integrable system (though by no means the 
only way!)
is in terms of a family $C\to B$ of spectral curves $C_b$.  In 
the simplest
situation, the abelian varieties $X_b$ are just the Jacobians 
$J (C_b)$.
More generally, $X_b$ may be realized as the Prym of $C_b$ with 
respect to an
involution, or as the generalized Prym with respect to a 
correspondence, cf.
\cite{D1}.
An example of such a system, due to Mukai \cite{Mu}, is obtained 
from the Jacobians of
a complete linear system $C\to B$ of curves on a symplectic surface
$(S,\barsig)$.  The symplectic form $\sigma$ on $X$ then corresponds 
via
Abel-Jacobi to the pullback to $C$ of $\barsig$.  Mukai studied the 
case where $S$ is a
$K3$ surface.  Taking $S = T^*E$ to be the cotangent bundle of some 
curve $E$
yields Hitchin's system on $E$.  One obtains many more examples by 
allowing $S$
to have a symplectic form $\barsig{_0}$ defined only on an open subset
$S_0\subset S$, and imposing appropriate restrictions on the 
intersection of the
spectral curves $C_b$ with $S-S_0$.  For example, we can take 
$S$ to be the total
space of the line bundle $\omega_E(D)$ for some effective divisor 
$D$ on $E$, to
recover Markman's system \cite{Mn} parametrizing meromorphic 
Higgs bundles on
$E$.  Note that this $S$ has a natural meromorphic 1-form 
$\barlam$, given away
from $D$ (i.e. where $S$ can be identified with $T^*E$) as 
the action 1-form
$\barlam = pdq$, whose differential is $d\barlam = \barsig{_0}$.  
It follows that
the ``tautological'' 1-form on $C$ (pullback of $\barlam$) 
induces the
Seiberg-Witten differential $\lambda$ on the integrable 
system $X =
J(C/B)\to B$.  The system discussed in \S 3 is of this form.

\medskip\noindent
{\bf Linearity: complexified Duistermaat-Heckman}

One feature of the Seiberg-Witten picture which we have not yet 
discussed is the 
linear dependence of the symplectic form on the particle 
masses.  From one point
of view, this linearity follows from existence of the 
Seiberg-Witten
differential:  the equation $\sigma = d\lambda$ implies that 
the cohomology
class of $\sigma$ is a fixed linear combination of the residues 
of $\lambda$,
which are just the masses.  On the other hand, linearity can 
also be interpreted
as arising from a complex analogue of the theorem of 
Duistermaat-Heckman
\cite{DH}.

Let $(X,\sigma)$ be a complex symplectic manifold (i.e. 
the symplectic form
$\sigma$ is holomorphic, of type (2,0)) equipped with a 
Hamiltonian action of a
complex reductive group $G$.  The Hamiltonian property of 
the action means that
there is a moment map $\mu: X\to \grg^*$ to the dual of 
the Lie algebra.
Ignoring bad loci, we will assume (or pretend) that there 
is a good
non-singular quotient $X/G$.  Let $G_a\subset G$ be the 
stabilizer of $a\in
\grg^*$ under the coadjoint action.  The symplectic 
quotient $X_a$ (or $X//_a G$)
is defined as
$$
X_a: = \mu^{-1}(a)/G_a = \mu^{-1}(\scrO_a)/G,
$$
where $\scrO_a$ is the (coadjoint) orbit of $a$.  It is 
again a complex
symplectic manifold, with form $\sigma_a$.  As $a$ ranges 
over orbits of a fixed
``type'' in $\grg^*$, the topology of $X_a$ remains locally 
constant, so 
cohomologies of nearby $X_a$'s of fixed type can be 
identified.  The
complex Duistermaat-Heckman theorem says that, under 
this identification,
the cohomology class $[\sigma_a]\in H^2(X_a)$ varies 
linearly with $a$.

Consider first the abelian case where $G = T$ is a torus $T =
Hom(\Lambda,\dbC^*)$.  The moment map then factors
$$
X\to X/T\to \bft^*,
$$
and we assume both maps are bundles, at least over some open 
$\bft^*_0
\subset \bft^*$.  The $T$-bundle $X\to X/T$ has first chern
class
$$
c_1\in H^2(X/T,\Lambda^*).
$$
For each $a$ in $\bft^*_0$, this restricts to a class 
(still denoted by $c_1$) in
$H^2(X_a,\Lambda^*)$.  The theorem in this case says that 
for $a, a_0\in\bft^*_0
= \Lambda\otimes\dbC$,
$$
[\sigma_a] -[\sigma_{a_0}] = \langle c_1, a-a_0\rangle.
$$
Because of possible monodromy of the fibers of $X/T$ over 
$\bft^*_0$, both sides
may be multivalued, so we should think of the linearity as 
being local in $a$.

Next, let $G$ be an arbitrary reductive group, and consider 
reduction at a
regular semisimple element $a$.  This means that the 
stabilizer $G_a$ is a
maximal torus $T\subset G$.  The coadjoint orbits of 
semisimple elements are
parametrized by the quotient $\bft^*/W$ of the dual 
Cartan by the Weyl group.
Let $\grg^*_0, \bft^*_0$ denote the open subsets of 
$\grg^*,\bft^*$ respectively
parametrizing regular semisimple elements.  The 
parameter space for regular
semisimple coadjoint orbits is then
$$
\bft^*_0/W \approx\grg^*_0/G.
$$
The moment map in this nonabelian situation no longer factors. 
The closest we
get to factorization is the commutative diagram:
$$
\spreadmatrixlines{1\jot}
\matrix
&& X_0 &&\\
&\swarrow &&\searrow\kern-4pt\raise4pt\hbox{${\ssize\mu}$} \\
\vspace{-4pt}
\hidewidth{X_0/G} &&&&\null\kern-5pt\grg_0^*\\
\vspace{-4pt}
&\searrow &&\swarrow &\\
&&&\hidewidth{\bold{t}_0^*/W = \grg_0^*/G} &
\endmatrix
$$
with $X_0: = \mu^{-1}(\grg^*_0)$.  Somewhat surprisingly, it is 
easy to reduce
this situation to the abelian case, by considering 
$Y: = \mu^{-1}(\bft^*)$.  The
symplectic form of $X$ restricts to a symplectic form on 
$Y$, and the $T$-action
on $Y$ is Hamiltonian with respect to this form, so there 
is a moment map
$\mu_T:Y\to \bft^*$ compatible with $\mu$ via the inclusions 
$Y\subset X$ and
$\bft^*\subset\grg^*$.  Now for $a\in\bft^*_0$ there is 
an obvious identification
of symplectic manifolds $X//_a G\approx Y//_a T$, so the 
dependence of the
symplectic class $[\sigma_a]$ of $X//_a G$ on $a \in t^*_0$ is 
again locally
linear.  (In the original case with $G$ compact, the 
multivaluedness is replaced
by discontinuities across walls, resulting in $[\sigma_a]$ being
linear on each Weyl chamber).

What happens when $a$ is allowed to degenerate?  If $G$ is compact,
$a$ remains semisimple (though not regular), and there is a
fairly obvious linearity for all such $a$ in a given open face
$\bold f_0$ of the Weyl chamber, i.e. a face $\bold f$ minus 
its proper
subfaces.  In the complex reductive case, there are more 
possibilities: the most natural orbit which corresponds
to an $a$ in a complexified open face $\bold f_0$ 
consists of elements which are not semisimple but regular, in
the sense that their centralizers have the minimal possible 
dimension. The closure of this orbit consists of 
a finite number of other orbits. The smallest of these orbits 
(the unique closed orbit, which is in the closure of each of 
the others) consists of the semisimple elements.
Although we will
not encounter non-regular elements in the integrable 
systems we need in the 
following section, it seems worthwhile to point out that 
linearity is indeed a
general phenomenon applicable to all strata.  

We identify $\grg, \bft$ with $\grg^*, \bft^*$ via an 
invariant form.
Consider an arbitrary stratum $S\subset \grg$ corresponding to 
a face $\bold f_0\subset\bold f\subset \bft$.  $S$ is a
bundle over (a finite quotient of the complexified) 
$\bold f_0$, and the fibers are (co)adjoint
orbits of a fixed type.  This bundle has local sections, 
which can be taken of
the form $n+\bold f_0$, where $n$ is some fixed nilpotent 
in the reductive
subalgebra $\bold l = \bold l_{\bold f}$ (= common 
centralizer of all 
elements of $\bold f_0$).  All elements $a = n+t$ of 
this section
have the same centralizer $\bold c$ in $\grg$.  
($\bold c$ is the centralizer in
$\bold l$ of $n$).  Since $\bold f$ is the center 
of $\bold l$, it is
contained in $\bold c$, which actually decomposes as 
$\bold c = \bold f + \bold
d$ with $\bold d$ semisimple.  Let $D\subset G$ be the 
subgroup corresponding to
$\bold d$.  The reduction to the abelian case now proceeds 
as follows:  We
observe that $\mu^{-1} (n+\bold f_0)/D$ is symplectic, and 
has a residual 
Hamiltonian action of the subtorus $F\subset T$ whose algebra 
is $\bold f$.  Let
$\mu_F$ be the moment map to $\bold f$ (which we identify with 
$\bold f^*$).
Then for $a = n+t$ as above, we have an isomorphism of symplectic
 manifolds:
$$
X//_a G\approx\mu_F^{-1}(t)/F = (\mu^{-1}(n+\bold f_0)/D)//_t F.
$$
Since the symplectic form on the right varies linearly 
with $t$ (which
parametrizes $a = n+t$), we conclude that $[\sigma_a]$ 
varies linearly
with $a$, as claimed.  

In a sense, there is linearity even as we degenerate
from one stratum to another.  We illustrate this phenomenon 
with the
following example.

The left action of $G$ on itself lifts to a Hamiltonian
action of $G$ 
on $X: = T^*G\approx G\times\grg^*$.  The quotient map 
$X\to X/G$ is
$$
\align
G \times\grg^*
&\to \grg^*\\
(g,a)
&\mapsto a,
\endalign
$$
while the moment map sends
$$
(g,a)\mapsto ad^*_g a.
$$
The symplectic quotients
$$
X//_a G = \scrO_a
$$
are the coadjoint orbits.   These fit, as discussed above, 
into a 
finite number of continuous families, each indexed by the 
$W$-quotient of some
open face.  Inside each family we expect linearity.  As we 
move between
families, the cohomology (and even the dimension) jumps.  
Yet there is a way
to express linearity uniformly for all $a\in\grg^*$.  Let 
$G_a,\grg_a$ 
be the stabilizers of $a$ in $G,\grg$ respectively.  Then $a$ 
restricts to
$\overline a\in \grg^*_a$ which is $G_a$-invariant, so it gives 
an element of
$(\grg_a^*)^{G_a}$, a bi-invariant $1$-formon $G_a$. This gives 
a class in $H^1(G_a,\dbC)$, which is then sent to 
$H^2(\scrO_a, \dbC)$ by
the transgression map (= the coboundary map 
$E^{01}_2\to E^{20}_2$ in the Leray
spectral sequence of the fibration $G\to\scrO_a$, with 
fiber $G_a$, where
$E^{pq}_2 = H^p(\scrO_a, H^q(G_a, \dbC))$).  The image 
is seen to be
$[\sigma_a]$, and the construction is obviously linear 
on each $\bold
f_0\subset\bold f\subset\bft$, since the stabilizers 
are actually constant there.
\bigskip
\bigskip

\head
\S3.  Which integrable system
\endhead

We have seen that all the structures on the 
quantum moduli space of 
$N=2$ SYM predicted by Seiberg and Witten can 
be realized in
terms of algebraically integrable systems.  The 
features of
``special geometry'' occur in any algebraically 
integrable system.  
The distinguished electric and magnetic charges 
appear as soon
as we fix a Seiberg-Witten differential, and 
perhaps the most 
natural way to get this is as tautological 
1-form on a family of curves
contained in the total space of a line bundle 
$\omega_E(D)$ of meromorphic 
differentials on some fixed curve $E$.  Since 
a quantum theory with 
$\beta = 0$ depends on a
parameter $\tau=\tau_{cl}$, it is natural to 
take $E$ to be the elliptic curve
corresponding to this $\tau$.  Finally, the 
linearity properties suggest that
our family of systems, parametrized by the 
masses, could arise as the family of
symplectic reductions of a larger 
symplectic manifold.

There remains the technicality of identifying 
a specific integrable system
for each type of theory, and testing that 
these theories have reasonable
physical properties.  In this section we 
first give several descriptions of
the system \cite{DW} which solves the $SU(n)$ 
theory with one adjoint, and
then describe some of the tests which confirm 
the choice.

\medskip\noindent
{\bf Meromorphic Higgs bundles}

Fix a curve $E$, an effective divisor $D$ on it, 
and an integer
$n\ge 2$.  The total space of the system is the 
moduli space $X$ of
equivalence classes of pairs $(V,\var)$, where $V$ 
is a rank-$n$ vector
bundle on $E$ with trivial determinant, 
$\varphi:V\to V\otimes \omega_E(D)$
is a traceless endomorphism with values 
in $\omega_E(D)$, and the pairs are 
subject to an appropriate semistability 
condition.  The base is 
$$
B: = \oplus^n_{i=2} H^0(E,(\omega_E(D))^{\otimes i}),
$$
and the map $\pi:X\to B$ sends a pair $(V,\varphi)$ 
to (the sequence of
coefficients of) the characteristic polynomial
$$
\det (t\cdot I-\varphi).
$$
(It is best to think of $t$ as a fiber coordinate 
taking values in
$\omega_E(D)$.)

The generic fibers of $\pi$ are abelian varieties.  
Consider the surface $S$
which is the total space of $\omega_E(D)$.  The 
points $b=(b_2,\dotsc, b_n)\in B$
parametrize a linear system of curves $C_b\subset S$, 
given by the equation
$t^n+b_2t^{n-2}+\cdots + b_n = 0$.  Let $\lambda$ 
be the tautological meromorphic
1-form on $S$.  Then the pair $(V,\varphi)$ 
determines a sheaf $L: = \text{ coker
}(\varphi-\lambda\cdot Id)$ on $S$ whose support is the 
curve $C_b$. In fact, for $(V,\varphi)$ 
generic, $L$ is a line bundle on
$C_b$.  Thus we have a map $\pi^{-1}(b)\to Pic(C_b)$. 
If we were working 
with $GL(n)$-bundles rather than $SL(n)$-bundles (remove 
the restrictions $\det V\approx
\scrO$, trace $(\varphi) = 0$), this would be 
an isomorphism of 
$\pi^{-1}(b)$ with a component of $Pic(C_b)$.  
The restrictions mean that we 
get instead a translate of Prym$(C_b/E)$.

This construction, modelled on Hitchin's 
system \cite{Hi}, works much more 
generally, cf. \cite{D2}.  In particular, 
the structure group $SL(n)$ can
be replaced by any reductive $G$.  The total 
space $X$ now parametrizes pairs
$(V,\varphi)$ with $V$ a $G$-bundle over $E$, 
and $\varphi$ a section of $ad(V)
\otimes\omega_E(D)$.  The base is replaced by
$$
B: = \oplus^r_{i= 1}H^0(E,(\omega_E(D))^{\otimes d_i}),
$$
where $r$ is the rank of $G$ and $d_i$ are the degrees 
of the basic
invariant polynomials on the Lie algebra $\grg$. 
 Once we choose
a representation $\rho$ of $G$, we get a family 
of spectral curves
$\{C_{\rho,b}\mid b\in B\}$ in $S$ as before.

As shown in \cite{D2}, \cite{Fa} and further 
references therein, the generic
fiber of $\pi$ can be identified as (a 
translate of) a certain distinguished
abelian variety Prym$(C_b)$ which is an 
isogeny summand in the Jacobian 
$J(C_{\rho,b})$ for each representation $\rho$.  
The most natural way to
view this is to abandon the spectral curves 
$C_{\rho,b}$ altogether, in favor
of the cameral cover $C_b\to E$, a Galois cover 
with group $W$, the Weyl
group of $G$.  This is the pullback to $E$, via 
the classifying map $E
\dashrightarrow \bold t/W\approx \grg/G$ determined 
locally in $E$ by the
point $b\in B$, of the $W$-cover 
$\bold t\to \bold t/W$.  ($\bold t$ is the
Cartan subalgebra, and the classifying map 
is really a map from the total space
of $(\omega_E(D))^*$ to $\bold t/W$.  The local 
form above depends on the choice
of a local section of the line bundle 
$(\omega_E(D))^*$, so is really defined only
modulo homotheties).  Now $J(C_b)$ decomposes 
according to 
the irreducible representations of the
Galois group $W$, and the distinguished Prym, 
Prym $(C_b)$, is the component
corresponding to the action of $W$ on the weights.  
(The question of the
dependence of the spectral Jacobian on $\rho$ was 
raised in \cite{AM} and 
studied, in the case of the Toda system, in \cite{Ml}. 
 A solution, under some
conditions, was given in \cite{Ka}.  
The case $E = \dbP^1$ was analyzed in
\cite{AHH} and \cite{B}, and the general case 
is in \cite{D2}, \cite{Me}).  In
this section we will use only the $SL(n)$ system, 
postponing the appearance of
other groups to \S 4.

The total space $X$ is not symplectic, but has a 
Poisson structure \cite{B}, \cite{Mn}.  In
general, a Poisson structure on an algebraic 
manifold $X$ is a Lie bracket
operation on $\scrO_X$ satisfying the properties 
of the usual Poisson bracket.
Equivalently, it can be given by a section $\alpha$ 
of $\wedge^2T_X$ satisfying
an integrability condition.  If $\alpha$ is 
non-degenerate on some open set
$X_0$, then its inverse $\sigma=\alpha^{-1}$ is 
a two-form on $X$, and
the integrability condition says that $\sigma$ 
is closed, i.e. $\sigma$
then gives a symplectic structure.  In general, 
a Poisson manifold
admits a natural foliation by submanifolds (whose 
conormal bundles consist of the
null spaces of $\alpha$) which inherit a symplectic 
structure.  A typical example
is $\grg^*$, which has the natural (Kirillov-Kostant) 
Poisson bracket, the
symplectic leaves being precisely the coadjoint orbits. 
More generally,
a Hamiltonian $G$-action on a symplectic manifold $Y$ 
induces a Poisson
structue on $X = Y/G$.  This is the ordinary quotient; 
the symplectic
quotients $X_a: = Y//_a G$, for $a\in\grg^*$, are 
recovered as the
symplectic leaves of $X$.  (The Kirillov-Kostant 
Poisson structure
on $\grg^*$ is recovered when we take $Y: = T^*G$).

Markman's construction of the Poisson structure is 
based on the
moduli space $\scrM_D$ of ``bundles with a level-$D$ 
structure'' on
$E: \scrM_D$ parametrizes equivalence classes of 
stable pairs $(V,\eta)$ 
where
$$
\eta:V_D\overset\sim\to{\longrightarrow}\scrO^n_D
$$
is a trivialization of the fibers of $V$ at 
points of $D$ (sending the given
volume element of $V_D$ to the standard one on 
$\scrO^n_D$, since our
structure group is $SL(n)$).  Consider the cotangent
 bundle $Y: = T^*\scrM_D$.  
By elementary deformation theory, its points are 
(equivalence classes of) 
triples $(V,\eta,\varphi)$ with $(V,\eta)\in\scrM_D$ 
and $(V,\varphi)\in X$,
a meromorphic Higgs bundle with poles on $D$.  The 
symmetry group $G$ is
the product of copies of $SL(n)$, one for each 
point of $D$.  (More precisely,
$G = \{g\in \text{ Aut}_D(\scrO^n_D)\mid\det g = 1\}$, 
allowing for possibly
non-reduced $D$.)  The $G$ action on $\scrM_D$ 
(with quotient $\scrM = \scrM_0$,
the moduli of stable bundles on $E$) lifts to a
 Hamiltonian action on $Y$:
$$
g:(V,\eta,\varphi)\mapsto (V, g \circ\eta,\varphi),
$$
with quotient $X = Y/G$, so we get a Poisson structure 
on $X$.  The moment map
$$
\mu:Y\to \grg^*\approx\grg
$$
sends $(V,\eta,\varphi)$ to $\varphi^\eta$, the 
element of End$_D(\scrO^n_D)$ 
which is conjugate via $\eta$ to the residue 
Res$_D\varphi\in\text{ End}_D(V_D)$.  
The symplectic leaf $X_a$, for $a\in\grg$, can 
therefore be identified as the
moduli space of pairs $(V,\varphi)\in X$ such that
 Res$_D\varphi$ is conjugate to the
given $a$.  In particular, each abelian variety 
$\pi^{-1}(b) \approx 
\text{ Prym}(C_b/E)$ is contained in one symplectic
 leaf $X_a$, and is a
Lagrangian subvariety there.

A number of special cases of this construction are 
well known.  When 
$D = 0$, we get Hitchin's system \cite{Hi} on
 $T^*\scrM$.  When $E$ has 
genus $0$, we get the polynomial matrices system 
studied in \cite{AHH} and
\cite{B}.  For more details and examples, 
see \cite{Mn} or \cite{DM1}.

For the SYM theory with adjoint matter, we take $E$ 
to be the elliptic
curve whose $\tau$ parameter (in the halfplane $\dbH$) 
is the coupling constant
of the theory.  We are looking for a 1-parameter 
family of integrable 
systems (parametrized by mass $m$), each having an 
$(n-1)$-dimensional
base and $(n-1)$-dimensional fibers.  One checks 
easily that $\dim(X_a)$ equals
the dimension of the adjoint orbit of $a$; so we 
are looking for orbits of the
smallest possible dimension $2(n-1)$.  There is 
essentially only one 1-parameter
family of possibilities:  $D$ consists of a single 
point $\infty\in E$, and $a$
is the conjugacy class of the diagonal matrix with 
eigenvalues $m(1,1,\dotsc, 1,
1-n)$.  This conjugacy class consists of all 
diagonalizable $SL(n)$ matrices $A$ with
$$\text{rank}(A-m\cdot Id)\le 1.$$
  (As $m\to 0$ there is one 
further, nilpotent, orbit
of the right dimension:  the conjugacy class of the 
elementary matrix $e_{12}$.
But this turns out to give a trivial integrable 
system, i.e. the same as for $a =
0$:  all the spectral curves are completely reducible).

Our system can be described rather explicitly.  The moduli 
space of degree 0,
rank $n$ vector bundles on $E$ can be identified with the 
symmetric product
$S^nE$: each equivalence class contains a unique 
decomposable bundle $V =
\oplus^n_{i=1} L_i$, with 
$L_i = \scrO_E(q_i-\infty)\in \text{ Pic}^0 E$.  
The moduli of
$SL(n)$-bundles is the fiber over $0$ of the 
Abel-Jacobi map $S^nE\to E$.  It is
a projective space
$$
\{q_1 + \dots + q_n \vert \sum q_i = n\infty\}
= \dbP H^0(E,\scrO(n\infty))\approx\dbP^{n-1}.
$$
In terms of the decomposition $V = \oplus L_i$, any 
Higgs field $\varphi$ on $V$ 
can be written as a matrix with entries
$$
\varphi_{ij}\in\Gamma(L_j\otimes L_i^{-1}(\infty)).
$$
Let $A: = \text{ Res}_\infty\varphi$.  Its diagonal 
elements are all zeroes, and $A
-m\cdot Id$ has rank $\le 1$, so $A$ is conjugate, 
via a diagonal matrix, to $m$
times $e = \Sigma_{i\not= j}e_{ij}$.  This means that 
each $(V,\varphi)$ in our
system is equivalent to a pair $(V,\varphi)$ with 
Res$_\infty(\varphi_{ij}) = m$
for $i\not= j$.  For a given $V$, this condition 
uniquely determines the
$\varphi_{ij}$ (as long as $L_i\not= L_i$).  The 
remaining free parameters are the
coordinates $q_i$ of $V = \oplus \scrO(q_i-\infty)$ and the
diagonal entries $p_i = \varphi_{ii}$.  Away from the 
diagonals $q_i = q_j$, these
$(q_i, p_i)$ can be identified with the canonical 
coordinates on $T^*\scrM$.

\medskip\noindent
{\bf The spectral curves}

The family of spectral curves $C_b$, for the 
fundamental,
$n$-dimensional representation of $SL(n)$,
can be described as a linear system on the
compactification $\overline S = 
\dbP(\omega_E(\infty)\oplus\scrO_E)$ of the
surface $S$, or on various birationally equivalent surfaces.  
More
precisely, it forms an affine open subspace of this 
linear system, the complement
consisting of curves which contain the section at 
infinity, $Z_\infty: =
\overline S-S$, with some multiplicity $k>0$.  
By removing this multiple component, we recover the spectral
curves for the analogous problem with $n$ 
replaced by $n-k$.  (If $A$ is an
endomorphism of a vector space $V$, and $W$ is 
an invariant subspace, then
$$\text{ rank}_W(A-m\cdot Id)\le \text{ rank}_V(A-m\cdot Id),$$ 
so the conjugacy class restricts correctly.)

We can write down an explicit equation
$$
p_n(t) = t^n + \Sigma_{k=2}^n f_kt^{n-k} = 0
$$
for the spectral curves.  Here $t$ is a coordinate on the 
fibers of
$\omega_E(\infty)$, and $f_k\in\Gamma(\scrO_E(k\infty))$. 
(We identify
$\scrO_E\approx \omega_E$).  Once we have such a $p_n$ for 
each $n$, the previous
remark and the tracelessness condition imply that the 
general spectral 
curve has equation
$$
C_u: p_n+\Sigma^{n-2}_{\ell=0} u_\ell p_\ell = 0,
$$
with arbitrary constants $u_\ell$ which can therefore be 
interpreted as
coordinates on the quantum moduli-space, the base of 
the integrable system.

The condition on $p_n$ which specifies the right conjugacy 
class at infinity is
that when we substitute $t = t'+\xi^{-1}$, where $\xi$ 
is a local coordinate on
$E$ at $\infty$, we get a function 
of $t'$ and $\xi$ with at most first 
order poles in $\xi$.
Write the equation of $E$ in Weierstrass form,
$$
E: y^2 = x^3 + bx-c.
$$
The local coordinate $\xi$ determines a sequence 
$x_k, k=2,3,\ldots$, of rational
functions on $E$ which are regular on $E-\infty$ 
and such that near $\infty$,
$$
x_k-\xi^{-k} = O(\xi^{-1}).
$$
We may as well use the local coordinate $\xi = x^{-1/2}$, 
and then the $x_k$ can
be read off the Taylor expansions:
$$
y=\xi^{-3}(1+b\xi^4-c\xi^6)^{1/2} 
=\xi^{-3}(1+\frac b2 \xi^4-\frac c2\xi^6-\frac{b^2}8\xi^8 + 
\frac{bc}4\xi^{10}-\frac{2c^2-b^3}{16}\xi^{12}+\ldots)
$$
and
$$
y^{-1}=\xi^3(1+b\xi^4-c\xi^6)^{-1/2} 
=\xi^3(1-\frac b2 \xi^4+\frac c2\xi^6+\frac{3b^2}8\xi^8 - 
\frac{3bc}4\xi^{10}+\frac{6c^2-5b^3}{16}\xi^{12}+\ldots).
$$
Explicitly, we get:
$$
\align
x_2 &= x\\
x_3 &= y\\
x_4 &= x^2\\
x_5 &= xy\\
x_6 &= x^3\\
x_7 &= x^2y-\frac b2 y\\
x_8 &= x^4\\
x_9 &= x^3y-\frac b2 xy + \frac c2y\\
x_{10} &= x^5\\
x_{11} &= x^4y-\frac b2 x^2 y+\frac c2 xy + \frac{3b^2}8 y\\
x_{12} &= x^6,
\endalign
$$
and so on. In terms of this basis, one checks by an easy induction
argument that
$$
p_n = t^n-\Sigma_{k=2}^n (-1)^k
(k-1)
\pmatrix
n\\
k\endpmatrix
x_k t^{n-k},
$$
so the first few polynomials are:
$$
\align
&1\\
&t\\
&t^2-x\\
&t^3-3xt + 2y\\
&t^4-6xt^2 +8yt-3x^2\\
&t^5 -10 xt^3 + 20yt^2-15 x^2 t+4xy\\
&t^6 -15 xt^4 + 40 yt^3-45 x^2 t^2 + 24 xyt-5x^3\\
&t^7-21 xt^5 + 70 yt^4-105x^2 t^3 + 84 xyt^2
-35 x^3 t+6( x^2y-\frac b2 y).
\endalign
$$
This is a bit cleaner than the formulas in \cite{DW}, 
mostly because we
set the sum of the roots in the Weierstrass 
equation of $E$ to 0.  A
more-or-less equivalent solution was given already 
in \cite{T}, where the
generating function for the $p_n$ was written in 
terms of theta functions and
exponentials on $E$.  A slightly different (but 
equivalent) set of polynomials
was given in \cite{IM1}.

\medskip\noindent
{\bf Elliptic solitons}

It is curious to note that the same system arises in
several rather different contexts.  Here is a brief outline
of one of these, following the work of Treibich and Verdier
\cite{T}, \cite{TV}.

Finite dimensional solutions of the KP hierarchy (= infinite
dimensional integrable system) are well understood.  They
are parametrized by Krichever data, consisting of a curve $C$,
a point $p\in C$, a line bundle $L$, and local 
trivializations.  Such
data determine a point in the total space of 
the hierarchy with the property that
the infinite collection of commuting KP flows 
emanating from the point
sweep out only a finite dimensional space, 
isomorphic to $J(C)$.  The 
curve $C$, the point $p$ and the trivializations 
remain fixed, while $L$
moves freely in $J(C)$.

An elliptic soliton is a finite dimensional 
solution of KP in which the
first flow, $KP_1$, already closes up, i.e. 
evolves on an elliptic curve
inside the possibly larger Jacobian $J(C)$.  
The direction of $KP_1$,
for points coming from Krichever data, is the 
tangent vector at $p$ to the
Abel-Jacobi image of $C$ in $J(C)$.  So the condition is that 
$J(C)$ should contain an elliptic curve $E$ which is tangent
at $AJ(p)$ to the Abel-Jacobi image $AJ(C)$.

Start with an $n$-sheeted branched cover $\pi:C\to E$ and 
a point $p\in E$ with
local coordinate $\xi$.  Let $\pi^*(p) = q_0+ \cdots + q_{n-1}$, 
and let
$\xi_i: = \xi\circ\pi$ near $q_i$.  The tangent vector 
$t_i: = \partial/\partial \xi_i$
lives naturally in 
$\Gamma(\scrO_{q_i}(q_i))$.  
Its image
$$
(AJ)_* t_i\in T_{AJ(q_i)}J(C)
$$
can then be identified with the image of $t_i$ under the 
coboundary map
$\delta$ for the sequence
$$
0\to\scrO_C\to \scrO_C(q_i)\to\scrO_{q_i}(q_i)\to 0.
$$
Similarly, the tangent at $\pi^* p$ to the image $\pi^*E$ 
of $E$ in 
$J(C)$ is given by $\Sigma^{n-1}_{i=0}\delta (t_i) =
\delta(\Sigma^{n-1}_{i=0} t_i)$.  The tangency condition 
then says
that $\delta(\Sigma^{n-1}_{i=0} t_i)$ is proportional 
to $\delta(t_0)$, and
therefore equal (by the residue theorem) to 
$\delta(n t_0)$.  This
happens if and only if there is some 
$f\in\Gamma(\scrO(\pi^* p))$
with
$$
\text{Res}_{q_i}(f\pi^*dz) =
\left\{
\aligned
1\quad  \quad \quad&1\le i\le n-1\\
1-n\quad \quad \quad &i = 0,
\endaligned
\right.
$$
where $dz$ is a non-zero closed 1-form on $E$.  This is 
equivalent
to saying that $\pi:C\to E$ is one of our spectral curves, 
and $f\pi^*dz$
becomes the Seiberg-Witten differential.

\medskip\noindent
{\bf Tests}

What makes this description of SYM with adjoint matter 
convincing for a
physicist, is that it can be ``tested'' by working out 
various limits and
special cases, and that it yields results consistent 
with the ``predictions'' 
based on the expected behavior of the physical system.  
We outline some of these
tests.

\medskip\noindent
{\bf Consistency with Seiberg-Witten}

For $n=2$, with given $E$ and $m$, our system produces 
a 1-parameter family
of double covers $C_u\to E$.  Each $C_u$ has genus 2, 
and the
integrable system is the family of 1-dimensional Pryms.  
Write the
equation of $E$ as 
$y^2 = \Pi^3_{i=1}(x-e_i) = x^3 + bx-c$, with 
$\sum e_i = 0$ as before.  The cover $C_u$ is given 
explicitly in terms of our polynomials $p_i$ as:
$$
0 = p_2 + up_0 = t^2-x+u.
$$
Its Galois group over the $\dbP^1$ with coordinate 
$x$ is $\dbZ_2\times\dbZ_2$, 
with involutions 
$\alpha:y\mapsto -y,\ , \beta: t\mapsto -t$, 
and $\alpha\beta$.
The quotients are a rational curve 
(with coordinate $t$); the original $E$; 
and another elliptic curve $E_u$, respectively.
$$
\vbox{
\epsfxsize=4.00in\epsfbox{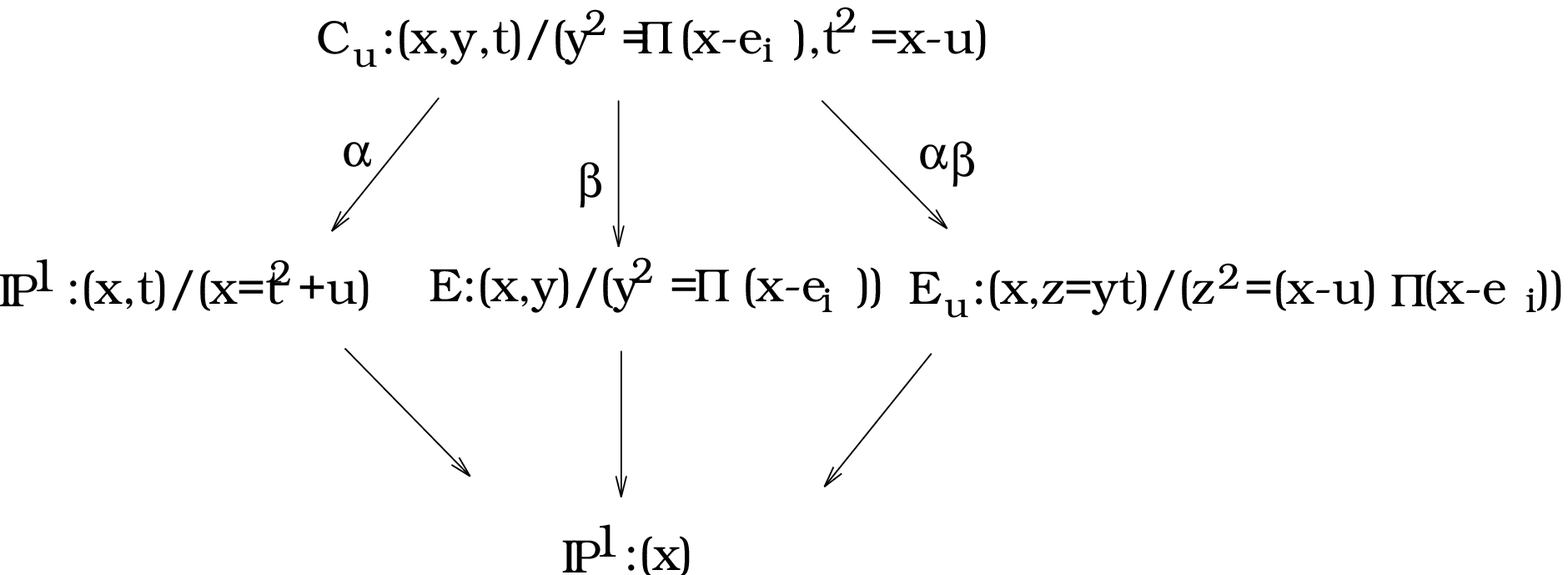}}
$$
\noindent 
The Prym can then be identified with $E_u$.

The original Seiberg-Witten solution is given directly as the
family of elliptic curves $E'_u$, double covers of $\dbP^1$ 
branched
at $e_i u + e_i^2 \quad (i=1,2,3)$ and 
$\infty$ (\cite{SW2}, eq. (16.17)).  These are
indeed isomorphic to our $E_u$: the transformation
$$
x\mapsto \frac{(u^2+b)x-c}{u-x}
$$
takes the $e_i$ and $u$ to $e_i u+e_i^2$ and
 $\infty$, respectively.

\medskip\noindent
{\bf Mass to zero: the N=4 limit}

The limit $m\to 0$ of $N=2$ SYM with gauge 
group $G$ and adjoint
matter is known to acquire a higher form of 
supersymmetry:  it is the 
$N=4$ SYM theory, with the same gauge group.  
(A physicist might rather
say that two of the four supersymmetries of $N=4$ SYM are 
broken when the theory is perturbed by the addition of a
mass term to the Lagrangian).
The extra symmetry of the $N=4$ theory means that 
``everything'' (i.e. the
ingredients of the special geometry package 
encountered in \S 2:  the quantum
moduli space, the Kahler metric, the prepotential...) 
can be computed
explicitly.  The outcome matches perfectly to the 
corresponding quantities for
Hitchin's system \cite{Hi}, i.e. the system of
 holomorphic (no poles)
Higgs bundles on $E$.  (This computation is in 
\S 2.3 of \cite{DW}).  So the
physics predicts that our system should specialize, 
as $m \to 0$, to Hitchin's for
$SU(n)$, which it of course does.  (This is not 
really an honest ``test'', since
we started out searching for deformations of 
Hitchin's system, so this
limit was built in from the beginning).

\medskip\noindent
{\bf $\tau$ to $\infty$: the flow to pure N=2}

One of the most interesting limits occurs as both $m$ and $\tau$ go to
infinity simultaneously, with $\tau$ proportional to $\log(m)$.
Equivalently, after an obvious rescaling, we can keep $m=1$ and let the
elliptic curve degenerate ($\tau \to \infty$) while adjusting the
coefficients $u_{\ell}$ in the spectral curves $C_u$.

The physics of this limit is easy to work out. Essentially, the mass
going to infinity means that in this limit the particle is so heavy
that it cannot be affected by the rest of the system; we get the pure
N=2 theory, with no mass. This theory was solved in [AF] and [KLTY].
The solution does not involve $\tau$, as we saw in the section on
adding matter: instead of $\tau$ having a limit (which was there called
$\tau_{cl}$), it now goes to $\infty$. The spectral curves found in [AF]
and [KLTY] have the explicit equation:  
$$ C'_b: w^2 = (z^n + b_2 z^{n-2} + \dots + b_n )^2 + 1 $$ 
in the $(z,w)$-plane, where the
coefficients $b_k$ play a role analogous to our $u_{\ell}$.

In [DW], we show that these curves are indeed recovered as the
appropriate limits of the spectral curves for the theory with adjoint
mass. This is done by substituting the $\tau \to \infty$ limit into the
explicit equations obtained above for the $C_u$, and performing the
relevant rescalings. Instead of repeating the derivation here, we will
describe the geometry of this degeneration.

The pure N=2 curve $C'_b$ is the fibre product over $\dbP^1$ (with
coordinate $v$) of the rational two-sheeted cover $w^2 = v^2 +1$ and
the (again rational, with coordinate $z$) n-sheeted cover $v=z^n + b_2
z^{n-2} + \dots + b_n$. In the adjoint theory, the analogue of the
double cover is just the elliptic curve $E \to \dbP^1$ (which becomes
rational in the limit $\tau \to \infty$), and $C_u$ is again an
n-sheeted branched cover of $E$, but it is no longer the fiber product
of $E$ with any $n$-sheeted cover of $\dbP^1$. The degree 2 map $\psi:
C \to \dbP^1_z$ and the degree n map $\dbP^1_z \to \dbP^1_v$ arise only
in the limit: $C$ becomes hyperelliptic.

The type of corank-1 degeneration involved in going from $C_u$ to
$C'_b$ can be seen more generally. Consider a degree n branched cover
$\pi : C \to E$ where $E$ is elliptic and $C$ has genus n. Up to
isogeny, the Jacobian decomposes $J(C) \approx P \times E$, where $P$
is the Prym. There are the Abel-Jacobi and Abel-Prym maps:

$$AJ: C \to J(C)$$
$$AP: C \to P$$

and their derivatives 

$$ \phi: C \to \dbP^{n-1}$$
$$ \psi: C \to \dbP^{n-2}.$$

\noindent Here $\phi$ is the canonical map and $\psi$ is the
composition of $\phi$ with linear projection from the point $0 \in
\dbP^{n-1}$  determined by the isogeny.

In the degeneration, we want each of $E,C$ to acquire a node: the
singular curves are $E_0, C_0$, with normalizations $E',C'$, and
singular points $p_0 \in E_0, q_0 \in C_0$ with branches $p_i, q_i,
i=1,2$ in the normalizations. The Jacobians $J(C), J(E)$ become
non-compact: $J(E_0)$ is a $C^*$, while $J(C_0)$ is a 
$C^*$-extension of  $J(C')$; but the Prym
$P$ remains compact. This type of degeneration occurs 
when $\pi' : C' \to E'$
 is totally ramified above each of the $p_i$: ${\pi'}^{-1}(p_i) = n q_i$.
 In this case we still have the Abel-Jacobi and Abel-Prym maps, going to the
 compactified Jacobians. The derivative $\phi$ is still the canonical
 map of the singular $C_0$; its image has a node at $0=\phi(q_0)$.
On the other hand, $\psi$ is no longer a morphism on $C_0$; its restriction $\psi':C' \to \dbP^{n-2}$ is now the canonical map of $C'$.

Within the $2n-2$-dimensional family of all covers $\pi : C \to E$ we restrict attention to the $n$-dimensional family of our spectral curves. Here there exist points $p \in E, q \in C$ such that $AJ(C)$ at $q$ is tangent to $\pi^*(E)$ at $p$. This is equivalent to saying that $AP: C \to P(C/E)$ is ramified at $q$, or that the linear system giving $\psi$ has a base point at $q$, or that $\phi(q)=0$.

Combining these two effects, we find that when the corank-1 degeneration described above occurs within the family of spectral curves, the two branches $\phi(q_1),\phi(q_2)$ are tangent at $0$. Equivalently, the canonical map $\psi'$ satisfies $\psi(q_1)=\psi(q_2)$; so $C'$ is hyperelliptic, $\psi$ is two-to-one, and ${\scr O}(p_1+p_2)$ is the hyperelliptic bundle. It follows immediately that the two maps, of degrees 2 and $n$ respectively, from $C'$ to $\dbP^{1}$ fit in a cartesian diagram. The corank-1 degeneration of the spectral curves (for the adjoint theory) therefore reproduces exactly the spectral curves of the pure $N=2$ theory.

\medskip\noindent
{\bf Higgs to $\infty$: symmetry breaking}

This time we fix both the curve $E$ and the mass $m$, and
instead we flow to the theory whose gauge group 
is some Levi subgroup
$G\subset SU(n)$.  Fix a diagonal matrix $\psi$, 
and consider the flow
sending a Higgs bundle $(V,\varphi)$ in our system 
to $(V,\varphi_s)$ where
$\varphi_s: = \varphi+s\psi$.  (As before we 
decompose $V = \oplus L_i$ and
let $\psi$ act diagonally.  The flow is 
well-defined only on the 
$n!$-sheeted cover of the system on which this 
decomposition exists).  The
physics predicts that the limit of the $SU(n)$
integrable system  as $s\to\infty$ along this 
flow should be the $G$-system, where $G\subset SU(n)$ 
is the 
centralizer of $\psi$.

For generic $\psi$, with distinct eigenvalues, this 
is clear:  After 
rescaling we have
$\varphi_s\sim\varphi'_s: = \psi + \frac1s\varphi$, 
so the 
limiting spectral curve is $\det(t-\psi) = 0$, which 
consists of $n$
copies of $E$ meeting at $\infty$.  The Prym is 
then the product
$J(E)^{(n-1)}$, independent of the initial 
$(V,\varphi)$.  This trivial
system corresponds to the Abelian structure 
group $U(1)^{n-1}$ which
is the maximal torus of $SU(n)$ and the 
centralizer of $\psi$.  

In general, 
$\psi$ will have, say, $k$ distinct 
eigenvalues with multiplicities $n_1,
\dotsc, n_k$, 
$\Sigma n_i = n$.  The 
centralizer is
$$
G = U(1)^{k-1}\times\Pi^k_{i=1} SU(n_i).
$$
Now the limit of the spectral curves of 
$\varphi'_s$ is non-reduced, consisting of 
$k$ components, the $i$-th with multiplicity 
$n_i$.  What we need is the
limit of the Jacobians.  This clearly has a 
factor (up to isogeny) from 
each component, and by rescaling around 
the $i$-th eigenvalue $\psi_i$ of
$\psi$ (i.e. setting $t-\psi_i = s^{-1}\tilde t$) 
we see that this 
factor is the Jacobian of the spectral curve
$$
\det(\tilde t-\varphi_i),
$$
where $\varphi_i$ is the $n_i\times n_i$ 
dimensional $i$-th block of $\varphi$.  
The claim is then
that this is one of the $SU(n_i)$ spectral 
curves, i.e. that the residue at
$\infty$ of the $i$-th block is 
conjugate to $m\cdot\text{ diag }(1,1,\dotsc, 1,1-n_i)$.  
An equivalent condition is: 
rank $(\varphi_i-m)\le 1$, and this
condition is clearly inherited by any 
block in $\varphi$, as required.

\medskip\noindent
{\bf Singularities}

Much of the effort in \cite{DW} was 
devoted to analyzing the
discriminant, the hypersurface 
$D\subset B$ over which the
spectral curves are singular.  Physically, 
that is where all the
interesting phenomena occur, and where many 
predictions can be
made.  For instance, as $\tau\to\infty$ 
(``at weak coupling''), the
limit of $D$ should be reducible:  this is a 
classical limit, where the
mass spectrum can be analyzed directly from 
the classical Lagrangian.
$D$ is the locus where some mass vanishes.  
This mass can belong
to either of the two types (vector or hyper) 
of $N = 2$ SUSY multiplets mentioned in section $1$ 
(``Adding matter''), 
and this gives rise to the two components.  In fact, one knows 
the explicit equations in $B$ 
for these components.  A Maple  assisted 
computation, using the explicit
equations obtained above for the spectral 
curves, confirmed the reducibility and
produced the precise predicted factors for $SU(3)$.  

For any $\tau$ (and any $n$), the most interesting 
locus in $B$ is the
$(n-1)$-fold self-intersection of $D$, a finite set.  
Physically, this
corresponds to the massive vacua, where the 
(low energy) gauge group
is broken to its maximal torus by monopole 
condensation.  A classification
result due to 't Hooft predicts that the possible 
massive vacua should correspond to
index $n$ subgroups of $\dbZ_n\times\dbZ_n$.

The corresponding result can be obtained directly 
from the integrable system.  
The general spectral curve has genus $n$, and we 
are imposing the maximal number,
$n-1$, of nodes.  The normalized spectral curve 
is then of genus 1, and
must be an unramified $n$-sheeted cover of the 
base $E$.  These covers are 
of course classified by the subgroups as above.  
Conversely, given any
$n$-sheeted cover $\pi:E'\to E$ and a point 
$\infty'$ above $\infty\in E$, there
is a unique meromorphic 1-form on $E'$ with 
residue $m$ at
$\pi^{-1}(\infty)\backslash\infty'$ and 
residue $(1-n)m$ at $\infty'$ (modulo 
the constant form $\pi^*dz$), hence a 
unique map of $E'$ to the total space
of $\omega_E(\infty)$ whose image is a 
spectral curve, as required.  (Changing
the choice of $\infty'$ amounts to a deck 
transformation, leaving the
spectral curve unchanged).

We also analyzed, assisted with some massive 
and miraculous-looking Maple 
computations, the locus of cusps in the $SU(3)$ 
theory, again confirming some
rather delicate predictions coming from the 
physics.  It seems that there is
some very beautiful geometry still waiting 
to be understood in the 
structure of the discriminant for 
the $SU(n)$ theories.

\head \S4.  Other integrable systems
\endhead

Having studied in some detail the system 
which solves one particular
family of SYM theories, we proceed to 
consider what is known for 
some of the other families.  Our 
discussion will include,

\roster
\item"{$\bullet$}"  The pure $N=2$ 
SYM for arbitrary gauge group 
$G$, and its relation to the Toda systems.

\item"{$\bullet$}"  The situation 
when fundamental matter is 
added to the pure theory.

\item"{$\bullet$}"  The theory with 
adjoint matter for arbitrary
$G$, and its relation to the elliptic 
Calogero-Moser system.

\item"{$\bullet$}"  A few words about 
the integrable systems arising in string theory.
\endroster

\medskip\noindent
{\bf Pure N=2 SYM}

The solutions by Seiberg and Witten of the pure $N=2$ SYM theory
\cite{SW1} and the theories with fundamental or adjoint matter
\cite{SW2} led rapidly to the creation of a whole literature of
solutions to specific SYM theories.  The general pattern was: 
postulate the general form of a family of curves which 
builds-in some
of the expected symmetries of the theory, then adjust 
coefficients to 
enforce all other available information:  further 
symmetries, limiting theories
which were already known, internal consistencies.  
All cases which were 
worked out this way, e.g. \cite{KLTY}, \cite{AF}, \cite{DS},
\cite{HO}, \cite{APS}, \cite{BL}, \cite{AS}, 
\cite{H}, \cite{AAG} and
further references therein, were based 
on postulating curves which are
hyperelliptic or closely related to hyperelliptics.

It was observed in \cite{GKMMM} that 
the curves obtained in
\cite{KLTY}, \cite{AF} for the pure 
$G = SU(n)$ theory are precisely the
spectral curves of the periodic Toda 
system.  The analogue for all semisimple
groups was found in \cite{MW}, which we follow here.

\medskip\noindent
{\bf The Toda system}

A simple way to describe the Toda system is 
via the construction of
Adler-Kostant-Symes.  Start with a 
decomposition $\grg = \bold a\oplus \bold b$
of a Lie algebra $\grg$ as vector space direct 
sum of two of its
subalgebras $\bold a,\bold b$.  This induces a 
dual decomposition 
$\grg^* = \bold a^*\oplus \bold b^*$, 
an injection 
$i: = \bold b^*\hookrightarrow \grg^*$, 
and a projection $\pi: \grg\twoheadrightarrow
\bold b$.  Let $\scrO_x\subset \bold b^*$ denote
the coadjoint orbit of some $x\in\bold b^*$. 
 The image $i(\scrO_x)\subset
\grg^*$ has, in general, nothing to do 
with the coadjoint orbit
$\scrO_{i(x)}$ of $i(x)$ in $\grg^*$.  
So $i^*$ of a $G$-invariant function
$f$ on $\grg^*$ can be a complicated, 
non-constant function on $\scrO_x$.  The 
AKS claim is that any two such, $i^*f_1$ 
and $i^*f_2$, Poisson commute with
respect to the natural (Kirillov-Kostant) 
symplectic structure 
$\sigma_{\scrO_x}$ on
$\scrO_x$.  Indeed, let $v_j\in\grg$ 
correspond to the differential
$df_j\vert_{i(x)}\in\grg^{**} = \grg$, 
so $i^*df_j$ is given by $\pi(v_j)$.
$G$-invariance of $f_j$ implies that $v_j$ 
is conormal to $\scrO_{i(x)}$, i.e.
that
$$
0 = \langle i(x),[\grg, v_j]\rangle.
$$
We set $a_j: = v_j-\pi(v_j)\in\bold a$, and compute:
$$
\align
\{ i^*f_1, i^*f_2\}_{\scrO_x}
&= \sigma_{\scrO_x}(\pi v_1, \pi v_2)\\
&= \langle x, [\pi(v_1), \pi(v_2)]_{\bold b}\rangle\\
&= \langle i(x), [\pi(v_1), \pi(v_2)]_{\grg}\rangle\\
&= \langle i(x), [v_1, \pi(v_2)]_{\grg}\rangle -
\langle i(x), [a_1, \pi(v_2)]_{\grg}\rangle\\
&= \langle i(x), [v_1, \pi(v_2)]\rangle
- \langle i(x), [a_1, v_2]\rangle 
+ \langle i(x), [a_1, a_2]\rangle.
\endalign
$$
We have just seen that the first two terms 
vanish when $f_1, f_2$ are
$G$-invariant, and the third term vanishes 
since it pairs $i(x)\in
\bold b^*= \bold a^\perp$ with 
$[a_1, a_2]\in\bold a$.  This proves AKS.

A shifted version of AKS is very 
convenient in applications.  Replace
the inclusion 
$i : \bold b^*\hookrightarrow \grg^*$ by its shifted
version $\epsilon + i$, where $\epsilon$ is 
an element of $\grg^*$
which annihilates $[\bold a, \bold a]$ and 
$[\bold b, \bold b]$.  The
above argument goes through verbatim.  
For example, take 
$\grg = {\bold s}{\bold l}(n)$,
$\bold a =$ strictly lower triangular matrices, 
$\bold b =$ upper triangulars,
and (after identifying $\grg$ with $\grg^*$, 
so that $\bold b^*$ becomes the
lower triangulars):
$$
\epsilon =
\pmatrix
0 &1\\
&0\\
&&\ddots\\
&&&0 &1\\
&&&&0
\endpmatrix,
x = 
\pmatrix
0\\
1 &0\\
&&\ddots\\
&&&0\\
&&&1 &0
\endpmatrix.
$$
The orbit $\scrO_x\subset \bold b^*$ 
consists of all matrices of the
form:
$$
\pmatrix
b_1 &&0\\
a_1 &\ddots\\
&\ddots\\
0 &&a_{n-1}b_n
\endpmatrix
$$
with $a_j\not= 0$ and $\Sigma b_j = 0$. 
The shifted orbit consists of tridiagonal matrices:
$$
\pmatrix
b_1 &1 &&&0\\
a_1 &\cdot &\cdot \\
    &\cdot &\cdot &\cdot\\
&&\cdot &\cdot &1\\
0 &&&a_{n-1} &b_n
\endpmatrix.
$$
This is a $2(n-1)$-dimensional symplectic 
manifold $(X,\sigma)$, isomorphic as
symplectic manifold to the cotangent 
bundle of the maximal torus $T\approx
(\dbC^*)^{n-1}$.  The restriction to 
it of the symmetric polynomials of $SL(n)$
(= coefficients of the characteristic 
polynomial) gives $n-1$ independent,
commuting Hamiltonians.  This is the 
classical Toda system.

An equivalent system is obtained if 
$\epsilon$ is replaced by $m\cdot\epsilon$
with $m\in\dbC^*$.  When $m = 0$, 
the Hamiltonians become the symmetric functions
of the $b_i$.  We get the cotangent 
bundle $X$ of $(\dbC^*)^{n-1}$, fibered by
its moment map over $(\dbC)^{n-1}$.  
Another variant is to take $\epsilon = 0$,
but replace $\bold a$ by the algebra 
$\bold s \bold o (n)$ of skew symmetric 
matrices.  Now
$\bold b^*$ sits in 
$\bold s \bold l(n)$ as the symmetric matrices, 
and the orbit consists of
symmetric tridiagonal matrices
$$
\pmatrix
b_1 &c_1\\
c_1 &b_2\\
&&\ddots\\
&&&b_{n-1} &c_{n-1}\\
&&&c_{n-1} &b_n
\endpmatrix.
$$
This is conjugate to the previous matrix 
with $a_i = c_i^2$, so we obtain a 
new version $(\tilde X,\tilde\sigma)$ of 
Toda which is just a $2^{n-1}$-sheeted
cover of the previous, basic version, $(X,\sigma)$.

The Toda system for a semisimple Lie algebra 
$\grg$ is obtained with $\bold b$ a
Borel subalgebra, $\bold a$ the nilpotent 
radical of the opposite Borel, $x =
\Sigma e_\alpha$ and $\epsilon = \Sigma e_{-\alpha}$, 
where the sums are over the
simple roots $\alpha$ with respect to $\bold b$, 
and the $e_{\pm\alpha}$ are
normalized root vectors.  The total space 
can again be identified with the
cotangent bundle of the maximal torus, 
a typical element
having the form 
$\op{\Sigma}\limits_{\alpha\text{ simple}}
(b_\alpha h_\alpha + a_\alpha e_\alpha + e_{-\alpha})$ 
with respect to a
standard basis $h_\alpha, e_{\pm \alpha}$.

\medskip\noindent
{\bf Periodic Toda}

The periodic Toda systems are obtained when 
the semisimple algebra $\grg$ is
replaced by a loop algebra $\grg^{(1)}$.  
This has an extended Dynkin diagram
consisting of the simple roots $\alpha$
 of $\grg$ plus the affine root
$\alpha_0$, corresponding to the highest 
root of $\grg$.  The general element of
the periodic Toda system can thus be written as
$$
\Sigma_\alpha (b_\alpha h_\alpha + a_\alpha e_\alpha + e_{-\alpha})
+ z e_{\alpha_0} + z^{-1}a_0 e_{-\alpha_0},
$$
where $z$ is the loop variable, and the product $a_0
\prod_\alpha a_\alpha = \mu$ is
a fixed non-zero constant.  
(For $\grg = \bold s \bold l (n)$, this is:
$$
\pmatrix
b_1 &1      &       &a_0z^{-1}\\
a_1 &\ddots &\ddots &          \\
    &\ddots &\ddots &1        \\
z   &       &a_{n-1}&b_n
\endpmatrix,
$$
with $\prod a_i = \mu, \Sigma b_j = 0$).

As in the case of Higgs bundles, the 
spectral curves $C_{\rho,b}$ depend on the
choice of a representation $\rho$ of 
$G$.  These Toda curves and their
relationships were first studied in 
\cite{Ml} and \cite{MS}, with some
clarifications added in \cite{MW1}.  
Each $C_{\rho,b}$ is a branched cover of the
$z$-line $\dbP^1_z$.  The latter plays 
the role for periodic Toda which the base
curve $E$ played for the meromorphic 
Higgs system in \S2.  The cameral cover in
this situation has a strikingly simple 
description.  The classifying map
$\dbP^1_z \dashrightarrow \bold t/W$, 
which is usually only locally defined, is now a
morphism for $z\not= 0,\infty$.  It 
factors through the double cover
$$
\align
\dbP^1_z
&\to \dbP^1_w\\
z &\mapsto w = z+\mu z^{-1}
\endalign
$$
which maps $\dbP^1_z \backslash
 \{ 0,\infty\}$ to $\dbA^1_w = \dbP^1_w\backslash
\{\infty\}$,
followed by an affine linear map
$$
\dbA^1_w\mapsto \bold t/W
$$
whose image in $\bold t/W$ is a straight line
 pointing in the direction of 
the invariant polynomial of highest degree, 
which is equal to the Coxeter number
$h_{\grg}$ of $\grg$.  More precisely, the 
algebra of polynomial functions on
$\bold t/W$ is isomorphic to the algebra of 
invariant polynomials on $\bold t$,
which is filtered by the degree.  The 
equations of the image of $\dbA^1_w$ are
obtained by fixing the values of all 
polynomials on $\bold t/W$ whose pullback to
$\bold t$ has degree $< h_{\grg}$.  
The induced $W$-cover of $\dbP^1_w$
then has $r = \text{ rank}(\grg)$ branch 
points in $\dbA^1_w$, and one at $w =
\infty$.  The choice of a Borel labels 
the finite branch points $w_\alpha$ by the
simple roots $\alpha$ in such a way that 
the monodromy at $w_\alpha$ is the
reflection corresponding to $\alpha$.  
The monodromy at $\infty$ is then the
product of these reflections, the Coxeter 
element.  The Cameral cover is the
pullback to $\dbP^1_z$.  It can be 
described as the $W$-cover obtained from the
trivial cover $W\times\dbP^1_z$ by 
making $r+1$ cuts labelled by the extended set
of roots $\{\alpha\}\cup \{\alpha_0\}$,
 each running from a point $z_\alpha$ to
$\mu z_\alpha^{-1}$, and pasting 
across each cut according to the corresponding
element of $W$.

In addition to the loop algebras $\grg^{(1)}$, 
there are twisted versions
$\grg^{(k)}$ where $k$ is the order of an 
automorphism of the Dynkin diagram of
$\grg$, and each of these also gives rise
 to a periodic Toda system.  The
basic observation of Martinec and Warner 
\cite{MW1} is that each of the known
families of pure $N=2$ SYM curves arises 
as spectral curves for one of these
twisted periodic Toda systems.  
The specific identification involves Langlands
duality (Olive-Montonen duality, 
to a physicist) exchanging the direction of
arrows in the extended diagram.  
The (untwisted) extended diagrams for types
$A,D,E$ are self-dual, yielding 
ordinary Toda systems.  The remaining cases yield
twisted algebras and systems:
$$
(B^{(1)}_r)^\vee = A^{(2)}_{2r-1},
(C^{(1)}_r)^\vee = D^{(2)}_{r+1},
(F^{(1)}_4)^\vee = E^{(2)}_6,
(G^{(1)}_2)^\vee = D^{(3)}_4.
$$
We will now describe explicitly the spectral 
curves for some of these systems.

\medskip\noindent
{\bf The spectral curves}

For the classical groups we can write the 
equation of the Toda spectral curves
$C_{F,u}$ for the fundamental representation $F$ as
$$
w = z + \mu z^{-1}, x^\epsilon w + p = 0
$$
where 
$$
p = p(x) = x^{n+1} + u_2 x^{n-1} + \cdots + u_{n+1}
$$ 
for type $A_n\,\,(G = SL(n+1))$, and 
$$
p = p(x^2) = x^{2n} + u_2 x^{2n-2_+}\ldots +
u_{2n}
$$ 
for types $B_n, C_n$, and $D_n$, and the $u_\ell$ 
are coordinates on the 
base.  Here $\epsilon = 0$ for $A_n, B_n, C_n$ 
and $\epsilon = 2$ for $D_n$.  In
general, $\deg_x(p)$ is the dimension of the 
representation, while the $w$ term
modifies the coefficient $u_{h_\grg}$ 
corresponding to the invariant of highest
degree.  For $A_n, B_n, C_n, D_n$, these 
Coxeter numbers are $n+1, 2n, 2n, 2n-2$.
(The $B_n$ case is slightly different 
since $0$ is a weight of $F$.  So the actual
``spectral'' curve has equation $x(w+p(x)) = 0$.  
We may safely ignore the $x$ factor, since
the distinguished Prym lives in the Jacobian of 
the other factor).

For $A_n$, the fundamental curve has genus $n$, 
so its Jacobian is the
distinguished Prym.  It can be described 
geometrically as the fiber product over
$\dbP^1_w$ of the double cover $w = z+\mu z^{-1}$, 
branched at $w = \pm
2\sqrt{\mu}$, and the degree $n+1$ cover $w = -p(x)$, 
with total ramification
over $w = \infty$ and $n$ simple branch points at 
the images of the roots of
$p'(x) = 0$.  The substitution 
$y = 2z + p(x)/2,\,\, 4\mu = \Lambda^{2(n+1)}$
converts it to
$$
y^2 = p^2(x)-\Lambda^{2(n+1)},
$$
the form obtained in \cite{AF} and \cite{KLTY}.  
(Here we have reinstated the
scale parameter $\Lambda$ of the asymptotically 
free theory.  Previously this was
set to 1).

The other spectral curves for $A_n$ are easily 
obtained from the fundamental one.
The curve for $\wedge^2 F$ is of degree 
$\frac{n(n+1)}2$ over the $z$ line, and
of genus $n(n-1)+1-[\frac{n+1}2]$.  The 
obvious correspondence between it and
$C_{F,u}$ produces a copy of $J(C_{F,u})$ 
in $J(C_{\wedge^2 F,u})$.  The $A_3$
case is described in \cite{Ml}, \cite{MS}:  
The tetragenal curve $C_{F,u}$ of
genus 3 gives rise via Recillas' trigonal 
construction (cf. \cite{D3}) to the
degree 6, genus 5 curve $C_{\wedge^2F,u}$, 
together with an involution whose
quotient is a trigonal curve $C'$ of genus 2. 
 The Jacobian of $C_{F,u}$ is
isogenous to Prym$(C_{\wedge^2 F, u}/C')$.  
(Actually, the degree 6 and 3 curves
produced by the trigonal construction have 
genera 7 and 4 respectively, but each
has 2 nodes, over $z = 0$ and $z = \infty$; 
the curves $C_{\wedge^2 F,u}$ and
$C'$, of genera 5 and 3, are their respective 
desingularizations).  A beautiful
picture of the $g = 11$ curve $C_{\wedge^2 F, u}$ 
for $A_4$ is given in
\cite{MW1}.  In contrast to these small curves, 
there is the 
gigantic cameral cover, of
degree $n!$ over $\dbP^1_z$ and of genus 
$1+\frac{n-1}{2(n+1)}(n+2)!$, whose
Jacobian contains those of all 
(irreducible components of) spectral covers.

For $D_n, C_{F,u}$ has equation $x^2 w+p(x^2) = 0$. 
 This depends on $n$
parameters, but has genus $2n-1$, too big.  
The quotient $C'$ by the involution
$x\mapsto -x$ is hyperelliptic of genus $n-1$; 
the $n$-dimensional piece is Prym
$(C_{F,u}/C')$, and it can be realized as 
the Jacobian of a curve $C''$, the
quotient of $C_{F,u}$ by the involution 
$x\mapsto -x, z\mapsto\mu z^{-1}$.  The
substitution
 $y: = x^2(z-\mu z^{-1}),\,\, 4\mu = \Lambda^{4(n-1)}$ 
takes the equation
of $C_{F, u}$ to 
$$
y^2 = p^2(x^2)-\Lambda^{4(n-1)} x^4,
$$
the form of the $SO(2n)$ SYM curve obtained in \cite{BL}.  
In terms of $s = xy, t
= x^2$, the quotient $C''$ becomes:
$$
s^2 = tp^2(t) -\Lambda^{4(n-1)} t^3.
$$

For type $C_n$, the fundamental Toda curve $C_{F,u}$ 
has equation $w + p(x^2) =
0, w = z+\mu z^{-1}$.  The substitution 
$y = 2z + p(x^2)$ changes this to
$$
y^2 = p^2(x^2) -4\mu.
$$
Just as in the $D_n$ case, $C_{F,u}$ has 
automorphism group $\dbZ_2\times\dbZ_2$
over the $t = x^2$ line:
$$
\vbox{
\epsfxsize=2.5in\epsfbox{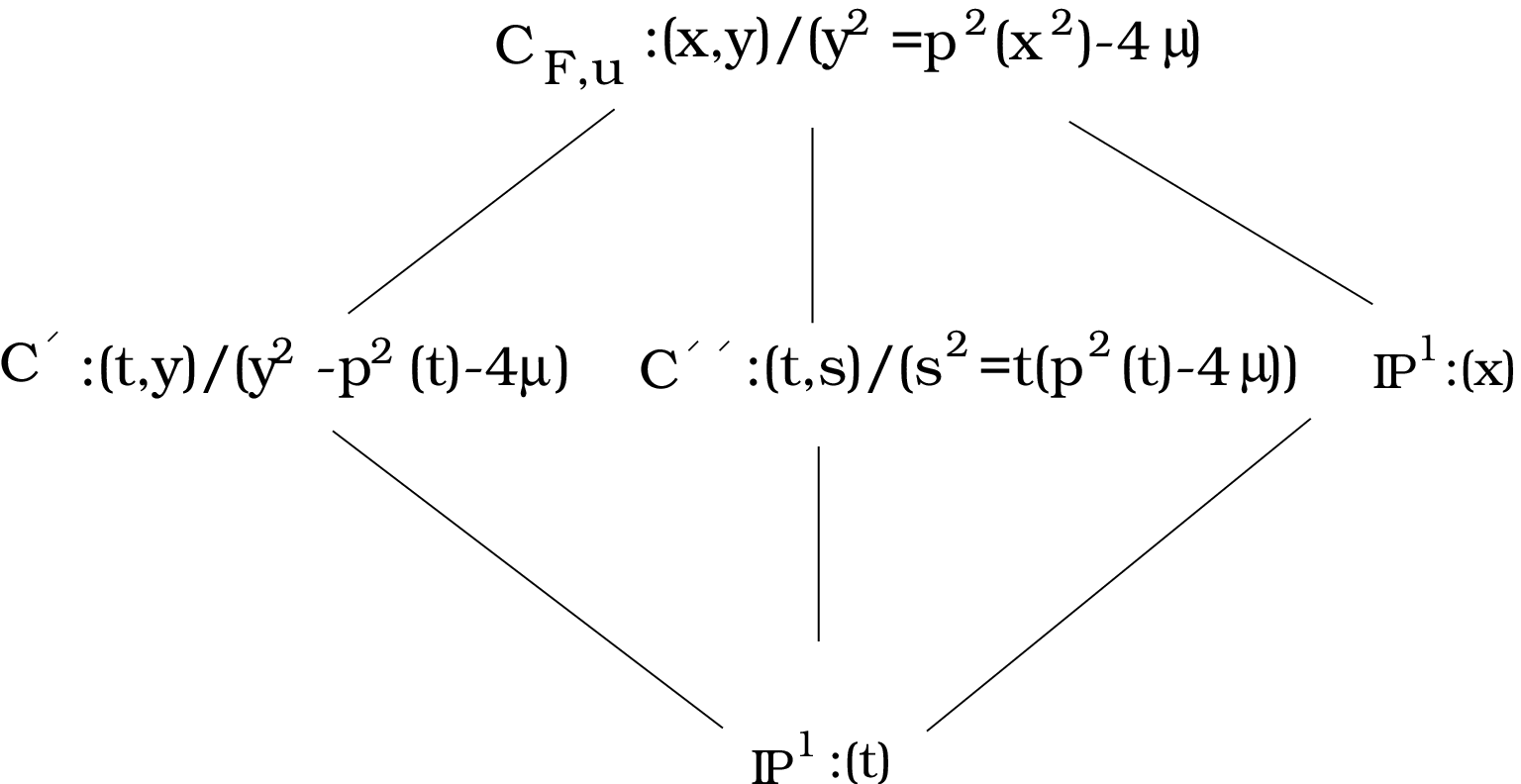}}
$$
\noindent and the $n$-dimensional piece is $J(C'')$.
But this is {\bf not} the right
abelian variety for the SYM theory.  Rather, we should 
consider the Toda curves
for $(C^{(1)}_n)^\vee = D^{(2)}_{n+1}$.  These have the form
$$
w^2 + p(x^2) = 0,\,\, w=z+\mu z^{-1}.
$$
The coordinates on the three $\dbZ_2$ quotients 
by the involutions
$x\mapsto\pm x, z\mapsto\pm \mu z^{-1}$ are now:
$$
\vbox{
\epsfxsize=1.5in\epsfbox{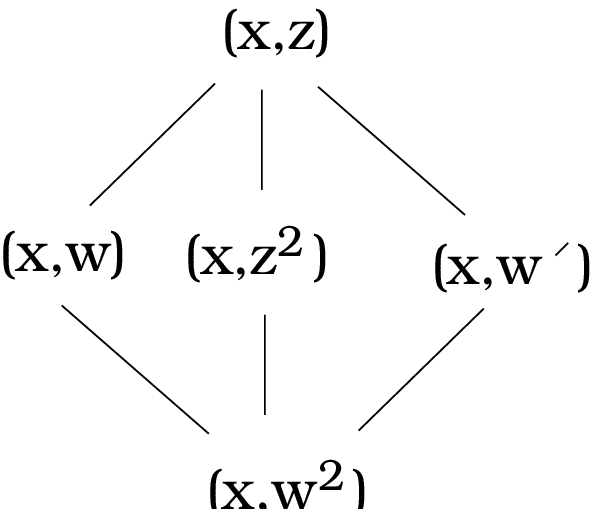}}
$$
\noindent 
where $w' = w-2\mu z^{-1} = z-\mu z^{-1}$. 
The corresponding genera are:
$$
\vbox{
\epsfxsize=1.5in\epsfbox{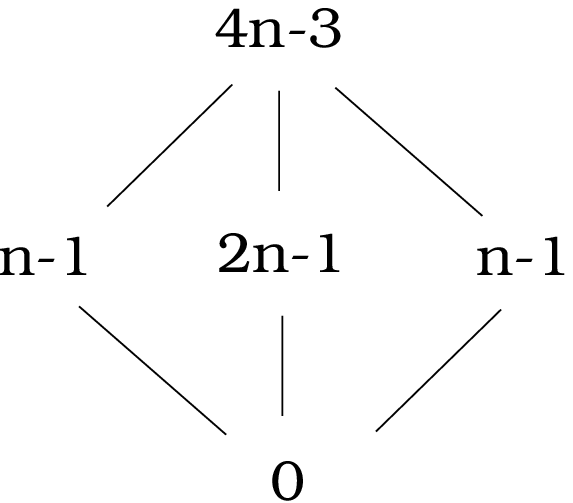}}
$$
\noindent
The middle curve has equation 
$v + \frac{\mu^2}v + p(x^2)-2\mu = 0$, where 
$v = z^2$.  It factors further:  in fact, 
it looks just like the $C_n$-Toda curve,
with $z,\mu, p$ replaced by $v,\mu^2, p-2\mu$ 
respectively.  So its Jacobian 
has an $n$-dimensional piece, 
which is the Jacobian of:
$$
s^2 = t((p(t)-2\mu)^2 - 4\mu^2) = t p(t) (p(t)-4\mu).
$$
This is essentially the $C_n$-curve obtained in \cite{AS}.

For $B_n$ the story is similar, but simpler.  
The dual Toda curve has equation
$x(z+\mu z^{-1}) + p(x^2) = 0$.  
The substitution $y = 2xz+p(x^2),\,\, 4\mu =
\Lambda^{2(2n-1)}$ converts this 
to the form obtained in \cite{DS},
$$
y^2 = p^2(x^2)-\Lambda^{2(2n-1)}x^2.
$$
Again this is of genus $2n-1$.  The genus $n$ quotient is
$$
s^2 = t(p^2(t) -\Lambda^{2(2n-1)}t).
$$

For $G_2$, we need Toda curves for
 $(G^{(1)}_2)^\vee = D^{(3)}_4$.  These are given
\cite{MW1} as:
$$
3(z-\mu z^{-1})^2 -x^8
+2u_2 x^6-[u^2_2+6(z + \mu z^{-1})] x^4 + 
[u_4 + 2u_2(z+\mu z^{-1})]x^2 = 0.
$$
The various quotients involve the functions:
$$
\align
w &= z+\mu z^{-1}\\
s &= x^2\\
u &= x(z-\mu z^{-1})\\
v &= x(w-s^2 + \frac{u_2}3 s).
\endalign
$$
Some quotient curves, indicated by 
the functions which generate their fields, are:
$$
\vbox{
\epsfxsize=1.5in\epsfbox{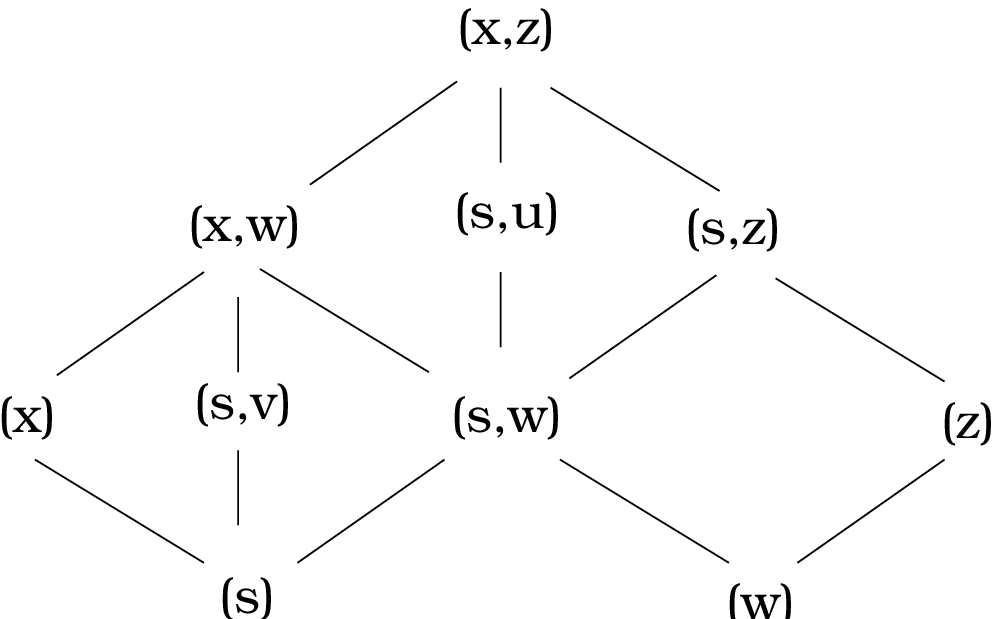}}
$$
The corresponding genera are:
$$
\vbox{
\epsfxsize=1.5in\epsfbox{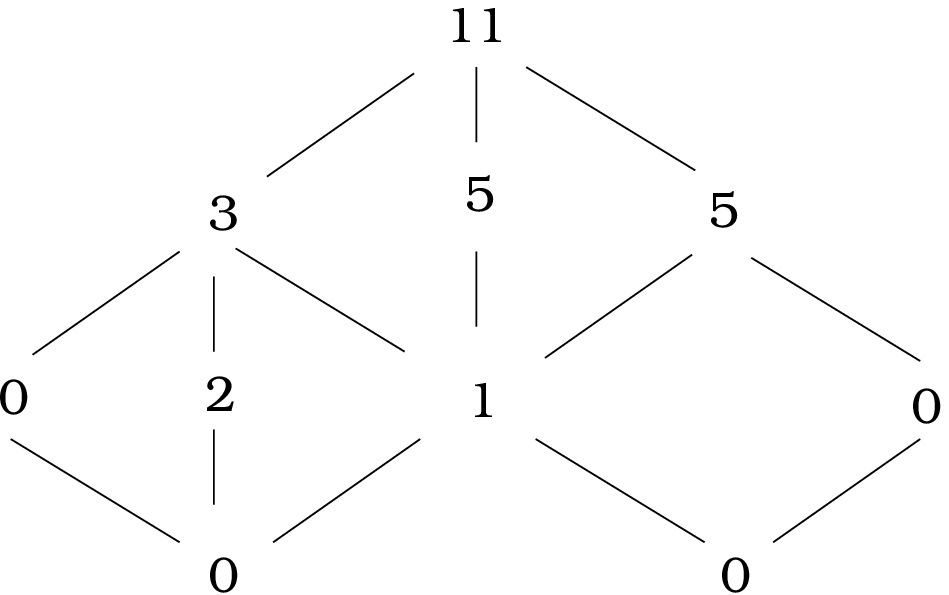}}
$$

I have not checked whether the genus 2 
curve with functions $(s,v)$ is the one 
corresponding to the $G_2$ theory.  An 
explicit family of genus 2 curves was
proposed in several works (see references 
in \cite{AAG}) for the $G_2$ theory,
but according to \cite{LPG} this family 
does not have good physical properties,
and does not match the Toda curves.

The $E_6$ curves were considered in 
\cite{LW}, and turn out to involve some
beautiful geometry.  The simplest curve 
is based on the 27-dimensional
representation of $E_6$.  The resulting 
27-sheeted cover of the $z$ line behaves
like the lines on a 1-parameter family of 
cubic surfaces.  Each of the 6 simple
roots corresponds to an ordinary double 
point acquired by the surface, so the
local monodromy at each finite branch 
point is the product of 6 disjoint
transpositions.  Lerche and Warner study 
this genus-34 spectral curve, and go on
to relate it to an integrable system coming 
from a string theory compactification on a
Calabi-Yau threefold which degenerates to a 
fibration of the $E_6$ ALE
singularity by a family of cubic surfaces.

\medskip\noindent
{\bf Fundamental matter}

The spectral curves for SYM with various 
gauge groups and numbers of particles
(``quarks'') in the fundamental representation 
have been determined by various
means.  A very incomplete list 
includes \cite{HO}, \cite{APS}, \cite{AS},
\cite{H}, \cite{MN}, \cite{KP}, \cite{DKP1},\cite{DKP2}. 
 Each of these solutions
is supported by substantial evidence, and yet 
the total picture seems far from
complete.  No really clear unifying principle 
is known, and there is no analogue
of the Toda integrable system which is 
responsible for the various curves.  There
is also a certain level of disagreement 
among comparable solutions, see the
discussion in \cite{MN} and \cite{DKP1}.  
(The curves involved in those solutions
are the same, but the parametrizations differ. 
An integrable system would of
course give a preferred parametrization).

\medskip\noindent
{\bf Adjoint matter}

As soon as \cite{MW} and \cite{DW} appeared, 
each with a class of integrable
systems which solves its respective 
version of SYM, the question arose whether
there was a common generalization:  
an integrable system for SYM with adjoint
matter and arbitrary gauge group $G$, 
going to the periodic Toda system in the
limit as $m$ and $\tau$ go together to $\infty$.

In \cite{Mc}, Martinec suggested that the 
elliptic Calogero-Moser system (cf.
\cite{OP}) may provide this common generalization. 
 There is such a system for
each semisimple $G$, and the case $G = SU(n)$ agrees
 exactly with the system in
\cite{DW}.  There is also a direct computation
 \cite{I1} showing that the $A_n$
and $D_n$ systems can degenerate to 
the respective Toda systems, and there is a
deformed family of integrable systems depending 
on some additional parameters
(the elliptic spin models, \cite{I2}). A way of obtaining 
the $A_n$ system by symplectic reduction is given in \cite{GN}.

One way to try to understand the geometry 
of these systems is as
follows.  Fix an elliptic curve $E$, 
a semisimple group $G$, and its maximal
torus $T$.  Let $\scrM, \scrM'$ denote 
the moduli spaces of semistable bundles
on $E$ with structue groups $G,T$ respectively. 
 A special feature of genus 1 is
that the natural map $\scrM'\to \scrM$ is  surjective, 
in fact it is a Galois
cover with group $W$.  (It is not true that 
the structure group of every
semistable $G$-bundle on $E$ can be reduced 
to $T$; but every $S$-equivalence
class of semistable bundles contains 
one whose structure group can be reduced).
Now let $X$ be the moduli space of Higgs 
bundles on $E$ with first order pole at
$\infty$, and let $X'$ be its finite cover 
induced by $\scrM'\to\scrM$.  (A
point of $X'$ is an equivalence class of 
pairs $(V,\varphi)$ where $V$ is a
$T$-bundle on $E$ and $\varphi$ 
is a meromorphic Higgs field on $V\otimes_T G$).
There is a ``moment map''
$$
\mu: X'\to\bold t^\perp/T,
$$
where $\bold t^\perp\subset \grg^*$ is the 
annihilator of the Cartan $\bold t$, and
the $T$-action is coadjoint.  Now {\bf all} 
the fibers of $\mu$ are symplectic,
of the minimal dimension $2r$ 
(where $r = \text{rank}(\grg)$), and at least for some of them
(but certainly not for all!), the restrictions of all the
$X$-Hamiltonians Poisson commute, 
giving a large family of small integrable
systems.  (Presumably, this involves 
an analogue of the AKS proof of integrability of
the Toda system).  As the elliptic curve 
degenerates to $\dbP^1$ with $0$ glued
to $\infty$, this construction specializes 
to the periodic Toda.

The question remains, though:  of all the 
possible $T$-orbits, which one
corresponds to SYM with adjoint matter?   
A necessary condition seems to be that
the $T$ orbit $\scrO$ needs to be $W$-invariant; 
otherwise the limit as $m\to 0$
will not be the Hitchin system $T^*\scrM$, 
but its finite cover $T^*\scrM'$.
(Or something in-between, if $\scrO$ is 
invariant under a subgroup of $W$.)  For
$SU(n)$, there is a unique such class 
(up to scalars), namely $\sum_{i\not=
j}e_{ij}$.  The natural analogue for 
arbitrary $\grg$ would be $\sum_\alpha
e_\alpha$, the sum over all roots.  
Unfortunately, this is ill defined:  the
$e_\alpha$ can be normalized only up to sign, and different signs give distinct
$T$ orbits.  For the orthogonal groups, 
there does not seem to be any
$W$-invariant $T$ orbit.  David Lowe 
has pointed out \cite{L} that for $Sp(n)$
there is an invariant orbit, namely 
the $T$ orbit of the defining matrix $J$.
(It is a signed sum of the long root vectors).
 This does give an integrable
system, but a rather boring one:  it is the $n$-fold 
product of the $SU(2)$
system, which does not flow to the pure
 $N=2\,\, C_n$ spectral curves described
earlier.

In sum, the elliptic Calogero-Moser systems 
(or their spin variants) seem to be a
promising general class of systems for 
realizing the SYM theory with adjoint
matter for arbitrary gauge group, 
but at the moment we lack a good geometric
understanding of these systems, 
e.g. via spectral curves, except in the case of
$SU(n)$.

\noi {\bf Stringy integrable systems?}

Quantum field theory is a limit of string 
theory compactified on a
Calabi-Yau threefold, and the SW duality 
can be realized as the limiting form of
various string dualities.  A series of 
recent papers, such as \cite{MW2},
\cite{KLMVW}, \cite{GHL1}, \cite{GHL2}, \cite{LW}, 
suggest that much of the QFT
structure studied in this paper - 
the spectral curves, integrable systems, SW
differentials and so on - arises as limit of 
corresponding stringy structures.  A
possible framework for studying the transition 
from stringy integrable systems to
the SW systems is proposed in \cite{F}, where 
SW-type integrable systems are
shown to arise from rigid special geometry, 
while the stringy type arises from
local special geometry.  From a purely 
mathematical point of view, it should be
interesting to understand the relations 
between the SW-type systems, whose basic
input is an elliptic curve; the Mukai 
system on a $K3$ surface, which degenerates
to Hitchin's on a curve (\cite{DEL}, \cite{DM1}); 
and the Calabi-Yau systems
constructed in \cite{DM2}.

\enddocument